\documentclass[11pt]{article}
\usepackage[letterpaper,margin=1in]{geometry}
\usepackage[T1]{fontenc}
\usepackage{lmodern}
\usepackage{amsmath,amssymb,mathtools,amsthm}
\usepackage{microtype,graphicx,booktabs,array}
\usepackage{algorithm,algorithmic}
\usepackage[algo2e,ruled,vlined,linesnumbered]{algorithm2e}
\SetKwInput{KwRequire}{Require}
\usepackage[round,authoryear]{natbib}
\usepackage{xcolor}
\definecolor{linkblue}{HTML}{244D73}
\usepackage[colorlinks=true,linkcolor=blue,citecolor=blue,urlcolor=blue]{hyperref}
\usepackage{pdflscape}
\usepackage{placeins}
\usepackage[font=small,labelfont=bf]{caption}
\newtheorem{theorem}{Theorem}[section]
\newtheorem{lemma}[theorem]{Lemma}
\newtheorem{proposition}[theorem]{Proposition}

\theoremstyle{remark}\newtheorem{remark}[theorem]{Remark}
\newtheorem*{remark*}{Remark}
\newcommand{\R}{\mathbb R}
\newcommand{\A}{\mathcal A}
\newcommand{\G}{\widehat{\mathcal G}}
\newcommand{\norm}[1]{\lVert#1\rVert}

\newcommand{\sgn}{\operatorname{sgn}}
\newcommand{\range}{\operatorname{range}}

\DeclareMathOperator*{\argmax}{arg\,max}

\graphicspath{{figures/}}
\hypersetup{pdftitle={Sign-Certified Broyden for Discounted Markov Decision Processes},pdfauthor={},pdfsubject={Preprint: global convergence and local superlinear convergence with nonunique optimal policies}}
\begin{document}
\title{A Globally Convergent Line-Search Quasi-Newton Method for Discounted Markov Decision Processes}
\author{
  Wangzhi Zhou, \ Yuanqiu Mo, \ and Wenwu Yu%
  \thanks{The authors are with the Department of Mathematics, Southeast University, Nanjing 211189, China (e-mails: \texttt{wangzhizhou474@gmail.com}, \texttt{yuanqiumo@seu.edu.cn}, \texttt{wwyu@seu.edu.cn}).}
}
\date{}
\maketitle
\begin{abstract}
In this paper, we develop a line-search quasi-Newton (LSQN) method for discounted Markov decision processes (MDPs). The method maintains an inverse model to approximate the inverse of the current policy Jacobian and refines it through an inner learning procedure that preserves existing secant conditions and achieves the prescribed accuracy in at most $n$ queries in exact arithmetic for an MDP with $n$ states, even when the initial inverse model is singular. A sign-based direction correction combined with an Armijo line search ensures monotonic decrease of the Bellman residual norm across policy changes without requiring a value-iteration fallback.
We prove global Q-linear residual convergence and local Q-superlinear convergence of the value iterates from arbitrary
finite initial value vectors and inverse models, without requiring
a unique optimal policy. Numerical experiments show competitive
performance against nine baselines, while public benchmarks
demonstrate that the inner learning procedure typically requires far fewer
than $n$ learning queries.
\end{abstract}

\section{Introduction}
Discounted Markov decision processes (MDPs) model sequential decision-making under uncertainty, whose objective is to maximize the expected cumulative discounted reward, with applications in inventory management, equipment maintenance, and communication systems \citep{Puterman1994}. Value iteration (VI) and policy iteration (PI) are two fundamental methods for finite MDPs, where VI repeatedly applies the Bellman operator and PI operates policy evaluation and improvement alternately \citep{Puterman1994,Bertsekas2012}. From the perspective of optimization,  VI is a fixed-point iteration,
whereas PI is a semismooth Newton method. Thus, this relation between the two schemes provides a direction for algorithm design
\citep{Gargiani2024SSN,Gargiani2024IPI}.

Based on the two basic methods, a lot of approaches are proposed to balance one-step cost and efficiency of iterates. Representative examples include accelerating VI by momentum \citep{GoyalGrandClement2022AVI}, anchoring \citep{LeeRyu2023AncVI}, and type-I Anderson \citep{ZhangODonoghueBoyd2020}, while another line reduces per-iteration cost of PI through approximate policy evaluation, such as modified policy iteration (MPI) performing a prescribed number of inner iterations before updating the policy \citep{PutermanShin1978,Bertsekas2012} and inexact policy iteration
\citep{Gargiani2023GMRES,Gargiani2024IPI} controlling the evaluation accuracy. Particularly relevant to our work is quasi-policy iteration (QPI) \citep{KolarijaniEsfahani2023QPI}, which employs structural constraints to approximate the policy Jacobian matrix. 

These works show different trade-offs in computational cost that cheap updates may lead to more iterations, while sufficient policy evaluation spends a lot. However, although strategies above enhance the efficiency of solving MDPs, how to balance computational cost with fast convergence remains an important question. Standard VI with cheap one-step costs may slow down with higher discount factors, while exact policy evaluation of PI would bring extra cost when policy changes continuously. MPI thus limits inner iterations to reduce cost, but may not match the difficulty of solving the current policy equation. For previous accelerated VI methods, they usually gain only a rate lower than the original one with global convergence from the first-order perspective. Second-order methods such as \cite{KamanchiDiddigiBhatnagar2022} face two challenges: smoothing errors and the high computational cost of solving a Newton system at each iteration. QPI reduces the computational cost by exploiting secant information, but does not guarantee sufficient accuracy in solving the policy equation. Inspired by the review of \citet[Section~5]{GrandClement2021Review}, we investigate quasi-Newton methods to accelerate VI, focusing on two key questions: \textit{How can policy evaluation accuracy be controlled to balance computational cost and convergence rate? And how can historical information be retained and reused to  improve the approximation of the current policy Jacobian inverse?} Quasi-Newton methods offer a natural approach by maintaining an inverse model to approximate the solution of the current policy equation. However, in MDPs, policy changes alter the underlying linear system, potentially making previously accumulated information less suitable for the current policy. Moreover, the behavior of the Bellman residual is not well controlled when iterates cross policy-region boundaries.

In this paper, we propose a line-search quasi-Newton (LSQN) method, serving as a scheme that sufficiently reuses information, and offering a middle ground between VI and PI. The algorithm keeps an inverse model across outer iterations, where an inner learning procedure is used to incorporate additional information about the current policy. Specifically, inner learning constructs secant conditions by the current policy Jacobian, while orthogonalization method \citep{BurdakovKamandi2018} is used to protect the secant conditions already satisfied, allowing new information to accumulate. Compared with exact PI, LSQN does not require solving the policy evaluation equation exactly in each iteration but keeps learning until a candidate step achieves the required precision, without requiring full recovery of the Jacobian inverse or a nonsingular initial inverse model, as required in traditional analyses \citep{BroydenDennisMore1973,BurdakovKamandi2018}. Unlike MPI, which prescribes the number of inner iterations, LSQN adaptively determines the amount of learning through a computable forcing residual test. The number of inner queries is uniformly bounded and can be much smaller when the current policy and retained inverse model provide a good approximation. 
Furthermore, to address the potential nonmonotonicity of the Bellman residual caused by policy changes, we construct a sign-based correction that ensures a common descent condition across near-active policy rows, together with an Armijo line search. In contrast to QPI, which employs a VI fallback or backtracking by embedding the value iteration residual directly into its search direction \citep{KolarijaniEsfahani2023QPI}, LSQN searches directly from the current value vector, utilizing a sign-based correction with an Armijo line search on a merit function, thereby guaranteeing global linear convergence without relying on VI updates. The secant-enhanced QPI method establishes local superlinear convergence via a Dennis--Mor\'e argument, without specifying a fixed order greater than one. By contrast, our forcing test yields an explicit local $Q$-order of at least $1+\nu$, where $\nu>0$ is the exponent in the forcing test.


Building on these techniques, we establish that LSQN achieves global Q-linear residual convergence and locally Q-superlinear value convergence from any finite initial value vector and inverse model, with explicit bounds on both convergence rates. The results establish faster local convergence than VI and avoid solving the policy-evaluation equation exactly. Compared with the global asymptotic convergence established for type-I Anderson acceleration \citep{ZhangODonoghueBoyd2020}, LSQN provides stronger convergence guarantees. Moreover, our analysis allows for multiple optimal policies.

Finally, we evaluate LSQN against nine baseline methods across four synthetic MDP families. The results show that LSQN requires fewer outer iterations than most competing methods. Despite its $O(n^3)$ worst-case computational complexity due to inner learning, LSQN achieves the second-best mean runtime rank, with a lower mean rank than the three QPI variants. Experiments on public benchmarks further show that the number of learning queries is typically far below the worst-case bound of $n$, keeping the practical computational cost low. The overall structure is outlined below.

Section~\ref{sec:method} introduces MDP preliminaries and the construction of the line search and the inner learning procedure. Sections~\ref{sec:global} and~\ref{sec:local} provide theoretical results of global linear convergence and the local superlinear convergence rate. Section~\ref{sec:learning} analyzes inverse-error reduction and the worst-case computational cost. Section~\ref{sec:experiments} reports the numerical experiments. The appendices give the envelope derivations and implementation
details.

\section{Line-Search Quasi-Newton Method}\label{sec:method}
In this section, we first introduce the preliminaries of MDPs. Next, we define the notion of a near-active set of actions, whose generalized-gradient envelope provides information about the step size required to guarantee descent. Finally, the proposed Line-search Quasi-Newton (LSQN) algorithm and the inner learning module are formally introduced. 
Throughout the paper, $\mathbb{R}$ denotes the set of real numbers. For $x\in\mathbb{R}^n$ and $A\in\mathbb{R}^{m\times n}$, we use $x_i$ (or $(x)_i$) and $A_i$ (or $(A)_i$) to denote the $i$-th component of $x$ and the $i$-th row of $A$, respectively. For $x \in \mathbb{R}^n$, $\operatorname{sgn}(x)$ denotes the sign vector obtained by applying the sign function componentwise, with $\operatorname{sgn}(0)=0$. For the Euclidean norm $\ell_2$, the infinity norm $\ell_{\infty}$, and cardinality, we denote by $\|\cdot\|_2$, $\|\cdot\|_{\infty}$, and $|\cdot|$, respectively.

\subsection{Bellman Residual and Merit Function}\label{sec:residual-facts}
Consider an MDP with states $i\in \mathcal{S} = \{1,2,\ldots,n\}$, a finite nonempty action set $\A$,
finite rewards $r_i^a$, and stochastic transition rows $P_i^a$.
The discount factor satisfies $0<\gamma<1$. Define
\begin{equation}
 q_i^a(v)=r_i^a+\gamma P_i^av,\qquad
 (Tv)_i=\max_{a\in\A}q_i^a(v)
 \label{eq:q_definitions}
\end{equation}
as the approximate $Q$-factor induced by choosing action $a$ in state $i$ and the maximum one among all actions, respectively, where $v \in \mathbb{R}^n$ denotes the state-value vector with $v_i$ estimating the expected discounted cumulative reward starting from state $i \in \mathcal{S}$, and $T$ is the Bellman operator. We then define the Bellman residual and the merit function as follows, respectively:
\begin{equation}
 F(v) = v- Tv, \quad \Phi(v)=\tfrac12\norm{F(v)}_2^2.
 \label{eq:merit}
\end{equation}
For a generic deterministic policy $\pi$, let $r_\pi$ and $P_\pi$ denote its reward vector and transition matrix, respectively,
and write $J_\pi=I-\gamma P_\pi$.
For a value vector $v$, we  denote $\pi_v$ as a greedy policy chosen by the rule $\pi_v(i)\in\arg\max_{a\in\mathcal A}q_i^a(v)$ for each state $i$, and consequently,
$F(v)=J_{\pi_v}v-r_{\pi_v}$. For the optimal value vector $v^*$, we define its induced optimal policy as $\pi^* := \pi_{v^*}$.
We also consider stationary randomized policies, represented by a matrix $\lambda=(\lambda_i^a)_{i,a} \in \mathbb{R}^{n\times|\mathcal A|}$, where $\lambda_i^a$ denotes the probability of choosing action $a$ in state $i$, satisfying
$\lambda_i^a\ge0$ and $\sum_{a \in \mathcal{A}}\lambda_i^a=1$ for all states $i \in \mathcal{S}$.
The corresponding transition matrix $P_\lambda$ has $i$-th row
$\sum_a\lambda_i^aP_i^a$, and we define
$J_\lambda=I-\gamma P_\lambda$.

\begin{lemma}
\label{lem:bellman-bounds}
Consider a finite discounted MDP with $0<\gamma<1$.
The Bellman operator $T$ defined in (\ref{eq:q_definitions}) is a $\gamma$-contraction in the
infinity norm and has a unique fixed point $v^*$ corresponding to the optimal value vector, i.e., $T(v^*) = v^*$. For every $v\in\mathbb R^n$,
\begin{equation}
 \|v-v^*\|_\infty
 \le \frac{\|F(v)\|_\infty}{1-\gamma},
 \label{eq:error}
\end{equation}
where F(v) is defined in (\ref{eq:merit}). Moreover, set $s_0=(1-\gamma)/\sqrt n$ and $L=1+\gamma\sqrt n$.
For every row-wise stochastic mixture $P_\lambda$,
the matrix $J_\lambda=I-\gamma P_\lambda$ satisfies
\begin{equation}
 s_0\|x\|_2
 \le \|J_\lambda x\|_2
 \le L\|x\|_2,
 \quad \forall x\in\mathbb R^n.
 \label{eq:matrix-bounds}
\end{equation}
In particular, $J_\lambda$ is invertible with
$\|J_\lambda\|_2\le L$ and
$\|J_\lambda^{-1}\|_2\le s_0^{-1}$.
\end{lemma}
The contraction property and the characterization of $v^*$ in \eqref{eq:error} follow from \citet{Puterman1994}, while \eqref{eq:matrix-bounds} holds trivially from $\|J_\lambda x\|_2\ge\|J_\lambda x\|_\infty
\ge(1-\gamma)\|x\|_\infty\ge s_0\|x\|_2$ and $\|P_\lambda\|_2\le\|P_\lambda\|_F\le\sqrt{n}$. This fundamental lemma converts Bellman-residual convergence
into value convergence by using the residual bound in \eqref{eq:error}, while the uniform matrix bounds as shown in \eqref{eq:matrix-bounds} govern the variations even when the selected policy changes. For further analysis, we provide the corresponding  bounds for Bellman residual differences under the $\ell_2$ and $\ell_{\infty}$ norms. Specifically, for all $x, y \in \mathbb{R}^n$,
 \begin{align}
 \norm{F(x)-F(y)}_\infty
 &\ge\norm{x-y}_\infty-\norm{Tx-Ty}_\infty\ge(1-\gamma)\norm{x-y}_\infty, \notag\\
 \norm{F(x)-F(y)}_2
 &=\left\|\int_0^1J_{\pi(t)}(x-y)\,dt\right\|_2\le L\norm{x-y}_2, \label{eq:two-bounds}
 \end{align}
where $\pi(t)$ in (\ref{eq:two-bounds}) selects greedy rows along $y+t(x-y)$ wherever the along-segment
derivative exists. 
The integral identity holds because $F$ restricted to the segment is a continuous, finite piecewise-affine map, which is absolutely continuous on each segment. Furthermore, in the case of persistent ties, the tied rows of $J_{\pi(t)}$ share the same derivative along that segment.

\subsection{Near-Active Descent Certificates}
Considering an MDP with an iterative method, the core objective is to reduce the Bellman residual until convergence, i.e., $F(v) \to 0$. In some special cases, perfect convergence properties can be achieved \citep{11535449} while ensuring the consistent decrease of the Bellman residual. However, without specific settings, iterative methods such as PI variants cannot guarantee the descent of the merit function $\Phi$ defined by (\ref{eq:merit}), even with exact policy evaluation.
To this end, we introduce the notion of near-active set and its generalized-gradient envelope to ensure descent across policy boundaries. Moreover, the upper and lower pairings are also defined to test whether a proposed search direction
satisfies the descent condition, and if not, to determine a sufficient correction.

At a non-optimal $v$, we define the near-active set as 
\begin{equation}
 \A_i^\varepsilon(v)=
 \{a\in\A:(Tv)_i-q_i^a(v)\le\varepsilon(v)\},
 \label{eq:width}
\end{equation}
where $\varepsilon(v)=\tau \|F(v)\|_2/\sqrt n$ given $\tau>0$.
Let $\Delta(\A_i^\varepsilon(v))$ denote the probability simplex on this set.
The corresponding generalized-gradient envelope is
\begin{equation}
 \G_\varepsilon(v)=
 \{J_\lambda^\top F(v):\lambda_i\in\Delta(\A_i^\varepsilon(v)),\
                         i=1,\ldots,n\}.
 \label{eq:envelope}
\end{equation}
Clarke calculus gives $\partial_C\Phi(v)\subseteq\G_\varepsilon(v)$.
Appendix~\ref{app:clarke} derives this inclusion from the finite-minimum
representation of each residual component. Moreover, based on Appendix~\ref{app:clarke}, for a candidate direction $d$, define its upper pairing as
\begin{equation}
 \begin{aligned}
 D(v;d)
:=\max_{h\in\G_\varepsilon(v)}h^\top d=\sum_{i \in \mathcal{S}}\max_{a\in\A_i^\varepsilon(v)}F_i(v)(d_i-\gamma P_i^ad),
 \end{aligned}
 \label{eq:D}
\end{equation}
where $F_i(v)$ denotes the $i$-th component of $F(v)$, and correspondingly the lower pairing is defined as 
\begin{equation}
 \begin{aligned}
 M(v;d)
:=\min_{h\in\G_\varepsilon(v)}h^\top d=\sum_{i \in \mathcal{S}}\min_{a\in\A_i^\varepsilon(v)}F_i(v)(d_i-\gamma P_i^ad).
 \end{aligned}
 \label{eq:M}
\end{equation}
The row-wise formulas follow from the independent choice of mixture
weights across states. Since $F_i(v)$ can have either sign, it must remain
inside each maximum or minimum.

\begin{remark}
The proposed two pairings serve different purposes. Given a direction $d$, the upper pairing $D(v;d)$ in \eqref{eq:D} tests such a direction against every vector in the envelope mentioned in \eqref{eq:envelope}, and $D(v;d)<0$ indicates descent even at policy boundaries. Instead, the lower pairing in \eqref{eq:M} provides a way to enforce this margin when a quasi-Newton proposal subsequently defined in Section \ref{sec:parameters} fails the test. At a non-terminal iterate,
let $b=\operatorname{sgn}(F(v))$. Lemma~\ref{lem:sign} subsequently 
shows that
\[
 M(v;b)\ge (1-\gamma)\|F(v)\|_1>0.
\]
Consequently, for any proposal $p$ and any $t\ge0$, we have
\[
\begin{aligned}
 D(v;p-tb) =\max_{h\in\G_\varepsilon(v)}
       \bigl(h^\top p-t h^\top b\bigr)\le D(v;p)-tM(v;b).
\end{aligned}
\]
Therefore, subtracting a suitable multiple of $b$ enforces the required
descent margin, as specified in \eqref{eq:correction}.
In summary, the upper pairing checks the proposal while the lower pairing fixes the direction.
\end{remark}

\subsection{Direction Construction and Line Search}\label{sec:parameters}
In the proposed LSQN algorithm, we use $k$ to index the outer iterations. At outer iteration $k$, an \textit{inverse model} is a matrix $H_k\in\mathbb R^{n\times n}$ that approximates the action of $J_k^{-1}=(I-\gamma P_{\pi_k})^{-1}$. Accordingly, the proposal $-H_kF(v_k)$ approximates the exact policy correction $-J_k^{-1}F(v_k)$. Before taking an outer step, we refine $H_k$ through an \emph{inner learning} module, detailed in Section~\ref{sec:completion} and Algorithm~\ref{alg:inner}, while keeping $v_k$, $\pi_k$, and $F(v_k)$ fixed. The inner learning procedure returns a modified matrix $\widehat H_k$ such that the resulting proposal $p_k=-\widehat H_kF(v_k)$ passes the forcing test in (\ref{eq:forcing}). At each outer iteration, if $\norm{F(v_k)}_2 = 0$, we terminate the algorithm since the optimal value vector has been obtained.

As will be shown in Lemma \ref{lem:completion}, the proposal
$p_k=-\widehat H_kF(v_k)$ with $\widehat H_k$ returned by the inner learning procedure will satisfy
\begin{equation}
 \norm{F(v_k) + J_kp_k}_2\le\eta_k \|F(v_k)\|_2,\qquad
 \eta_k=\min\{\bar\eta,c\|F(v_k)\|_2^\nu\},
 \label{eq:forcing}
\end{equation}
where $0<\bar\eta<1$, $c>0$, and $\nu>0$.
Specifically, $\bar\eta$ controls global progress, while the vanishing term $c\|F(v_k)\|_2^\nu$
sets the local rate.

Let $\theta \in (0, 1-\bar\eta]$, $\tau \in (0, 1-\theta)$ and $b_k=\operatorname{sgn}(F(v_k))$. $d_k$, the direction of line search in the $k$-th iteration is constructed based on the proposal such that
\begin{equation}
 d_k=
 \begin{cases}
 p_k,
 & \text{if } D(v_k;p_k)\le -\theta \|F(v_k)\|_2^2,\\[6pt]
 p_k-\displaystyle
 \frac{D(v_k;p_k)+\theta \|F(v_k)\|_2^2}{M(v_k;b_k)}\,b_k,
 & \text{otherwise}.
 \end{cases}
 \label{eq:correction}
\end{equation}
In particular, the proposal is used as the trial direction when it satisfies the descent certificate. Otherwise, a sign correction along $-b_k$ is applied to enforce the required margin. 
Lemma~\ref{lem:sign} establishes that $M(v_k;b_k)>0$
at every non-terminal iterate and that
$D(v_k;d_k)\le-\theta \|F(v_k)\|_2^2$ in both cases. However, although the directional certificate controls the first-order change in the merit function, it does not guarantee sufficient decrease over a full step. We therefore initialize $\alpha = 1$ and successively reduce it by a factor of $\beta\in(0,1)$, i.e., $\alpha = 1,\beta,\beta^2,\ldots$, until the following Armijo-type
sufficient decrease condition holds:
\begin{equation}
 \Phi(v_k+\alpha d_k)
 \le\Phi(v_k)+\sigma\alpha D(v_k;d_k),
 \label{eq:armijo}
\end{equation}
where $\sigma \in (0,1/2)$ denotes an Armijo coefficient that specifies the required fraction of the decrease predicted by the certificate. In short, the factor $\beta$ controls how quickly the trial step
is shortened, while $\sigma<1/2$ is used
later to establish eventual acceptance of full steps. Denote $\alpha_k$ as the accepted step size and $v_{k+1}=v_k+\alpha_kd_k$, with $s_k=v_{k+1}-v_k$ and $y_k=F(v_{k+1}) - F(v_k)$. At the end of each outer iteration, the inverse model is updated using the inverse Broyden update:
\begin{equation}
 H_{k+1}=\widehat H_k+
       \frac{(s_k-\widehat H_ky_k)y_k^\top}{y_k^\top y_k}.
 \label{eq:trajectory}
\end{equation}
Equation~\eqref{eq:two-bounds} makes the denominator in (\ref{eq:trajectory}) nonzero for every
nonzero accepted step. This update also enforces $H_{k+1}y_k=s_k$
\citep{Broyden1965}.

\setcounter{algocf}{\value{algorithm}}

\begin{algorithm2e}[htbp]
\DontPrintSemicolon
\caption{Line-Search Quasi-Newton for a discounted MDP}
\label{alg:outer}

\KwRequire{
Rewards $r_i^a$, transitions $P_i^a$, discount $0<\gamma<1$,
initial value vector $v_0$, initial inverse model $H_0$, parameters from Section~\ref{sec:parameters}.
}

\For{$k=0,1,\ldots$}{
    Compute $q_i^a=r_i^a+\gamma P_i^av_k$ and
    $(Tv_k)_i=\max_a q_i^a$.
    Select greedy policy as
    $\pi_k(i)\in\argmax_a q_i^a$
    and set $F(v_k) =v_k-Tv_k$.\;

    \If{$\|F(v_k)\|_2=0$}{
        \textbf{return} $v_k$ and $\pi_k$.\;
    }

    Compute $J_k=I-\gamma P_{\pi_k}$,
    $\A_i^\varepsilon(v_k)$ by \eqref{eq:width},
    and $\eta_k$ by \eqref{eq:forcing}.\;

    $(\widehat H_k,p_k)\leftarrow
    \textsc{Inner-Learning}(J_k,F(v_k),H_k,\eta_k)$,
    Algorithm~\ref{alg:inner}.\;

    \If{$D(v_k;p_k)>-\theta \|F(v_k)\|_2^2$}{
        Compute $b_k=\sgn(F(v_k))$ and $d_k = p_k-\displaystyle
 \frac{D(v_k;p_k)+\theta \|F(v_k)\|_2^2}{M(v_k;b_k)}\,b_k$.\;
    }
    \Else
    {    Set $d_k=p_k$.\;
    }
    Set $\alpha=1$.\;


    \While{
        $\tfrac12\norm{F(v_k+\alpha d_k)}_2^2
        >\tfrac12\norm{F(v_k)}_2^2+\sigma\alpha D(v_k;d_k)$
    }{
        Set $\alpha\leftarrow\beta\alpha$.\;
    }

    Set $v_{k+1}=v_k+\alpha d_k$,
    $s_k=\alpha d_k$,
    and $y_k= F(v_{k+1}) - F(v_k)$.\;

    Compute $H_{k+1} = \widehat H_k+
       \frac{(s_k-\widehat H_ky_k)y_k^\top}{y_k^\top y_k}$.\;
}
\end{algorithm2e}

\setcounter{algorithm}{\value{algocf}}

Algorithm~\ref{alg:outer} evaluates $Tv_k$ to form the Bellman residual, rather than using it as a fallback iterate or trial-path anchor, as is commonly done in safeguard strategies to guarantee a linear convergence rate \citep{KolarijaniEsfahani2023QPI}. In contrast, the linear convergence rate established in this paper is ensured by the direction certificate and line search, which together enforce a uniform decrease of the merit function in \eqref{eq:merit}.
\FloatBarrier
\subsection{The Inner Learning Procedure}\label{sec:completion}
The inner learning procedure improves the proposal \(p_k\) by learning the exact inverse actions of the current policy Jacobian \(J_{k}\). In each step, we evaluate \(g = J_{k}u = u - \gamma P_{\pi_k}u\) for a chosen direction \(u\). We call this evaluation a \textit{learning query}. Each query appends a new column to the learned representation of the inverse model, with the corresponding image orthogonalized against the previously learned image space. This adds an independent inverse secant condition to the model without overwriting previously learned ones, providing the exact relation \(J_{k}^{-1}g = u\). This further leads to a finite-completion mechanism: once the learned images span a basis of \(\mathbb{R}^n\), the learned inverse model matches the true policy inverse over the entire space.

Given the outer iteration $k$ and use $j$ to index iterations in the inner learning procedure, initialize $H^{(0)}=H_k$, and start
with empty matrices $U^{(0)},Q^{(0)}$ and $S^{(0)}$.
After $j$ queries these matrices have $j$ columns. The columns of $U^{(j)} \in \mathbb{R}^{n \times j}$
are orthonormal primal directions, those of $Q^{(j)}\in\mathbb{R}^{n\times j}$ are the corresponding orthonormal images, and $S^{(j)}\in\mathbb{R}^{n\times j}$ stores their exact inverse images. The following relations are maintained:
\begin{equation}
 J_kS^{(j)}=Q^{(j)},\qquad H^{(j)}Q^{(j)}=S^{(j)}.
 \label{eq:invariant}
\end{equation}
The first relation identifies the stored inverse actions, since
$S^{(j)}=J^{-1}_kQ^{(j)}$. The second requires the updated model to reproduce
them. Thus, for every coefficient vector $a\in\R^j$,
\[
 H^{(j)}Q^{(j)}a=S^{(j)}a=J^{-1}_kQ^{(j)}a.
\]
In particular, when \(j=n\), the orthonormality of the columns of $Q^{(n)}$ gives $Q^{(n)}(Q^{(n)})^\top=I$, and hence
\[
 \begin{aligned}
 H^{(n)}
 =H^{(n)}Q^{(n)}(Q^{(n)})^\top=S^{(n)}(Q^{(n)})^\top=J^{-1}_kQ^{(n)}(Q^{(n)})^\top=J_k^{-1}.
 \end{aligned}
\]
The proposal $p^{(n)}=-H^{(n)}F(v_k)$ then satisfies
$F(v_k)+J_kp^{(n)}=F(v_k)-J_kJ_k^{-1}F(v_k) = 0$, so $H^{(n)}$ passes \eqref{eq:forcing}.
Note that full inverse recovery is only the worst-case endpoint of this argument. The inner learning procedure terminates as soon as the forcing condition is satisfied.

The construction below obtains the inverse-image pairs through products with $J_k$, without explicitly computing $J_k^{-1}$. As shown in Lemma~\ref{lem:completion}, orthogonalizing each new image against the previously learned images allows its secant condition to be imposed without altering the model on the previously learned image space. Moreover, the relation $\range(S^{(j)})=\range(U^{(j)})$ is preserved throughout the procedure.

Set $p^{(j)}=-H^{(j)}F(v_k)$. If $p^{(j)}$ satisfies the forcing test in \eqref{eq:forcing}, the inner learning procedure terminates. Otherwise, project $p^{(j)}$ onto the orthogonal complement of $U^{(j)}$ and normalize the projection if it is nonzero. If the projection is zero, instead use the normalized projection of a coordinate vector with the largest projected norm. In either case, the resulting $u_j$ is a unit vector orthogonal to $U^{(j)}$ and is appended to $U^{(j)}$.

Next, let $g_j=J_ku_j$.
The component of $g_j$ orthogonal to the previously learned image space is given by $g_j-Q^{(j)}(Q^{(j)})^\top g_j.$
Its normalized version, which is appended to $Q^{(j)}$, is
\begin{equation}
z_j
=\frac{(I-Q^{(j)}(Q^{(j)})^\top)g_j}
{\norm{g_j-Q^{(j)}(Q^{(j)})^\top g_j}_2}.
\label{eq:orthogonalize}
\end{equation}
From (\ref{eq:invariant}), the corresponding inverse image $w_j$ appended to $S^{(j)}$ is given by $J_kw_j=z_j$, i.e., 
\begin{align}
 w_j= J_k^{-1}z_j &= J_k^{-1} (g_j-Q^{(j)}(Q^{(j)})^\top g_j)/\norm{g_j-Q^{(j)}(Q^{(j)})^\top g_j}_2 \notag\\
&= J_k^{-1}(J_ku_j - J_kS^{(j)}(Q^{(j)})^\top g_j)/ \norm{g_j-Q^{(j)}(Q^{(j)})^\top g_j}_2
 \notag\\
 &= (u_j-S^{(j)}(Q^{(j)})^\top g_j)/\norm{g_j-Q^{(j)}(Q^{(j)})^\top g_j}_2.
 \label{eq:inverse-image}
\end{align}
Together with \eqref{eq:orthogonalize} and \eqref{eq:inverse-image}, the inverse model is refined in the inner learning procedure as
\begin{equation}
 H^{(j+1)}=H^{(j)}+
       \bigl(w_j-H^{(j)}z_j\bigr)z_j^\top.
 \label{eq:inner-update}
\end{equation}
Since the newly appended image is orthogonal to all previously learned images, the rank-one update in \eqref{eq:inner-update} enforces the new secant condition without altering the previously satisfied ones. This preservation mechanism follows the classical projected-secant
principle \citep{BurdakovKamandi2018}. In particular, $H^{(j+1)}z_j
=H^{(j)}z_j+\bigl(w_j-H^{(j)}z_j\bigr)z_j^\top z_j
=w_j,$
whereas, for any previously learned image $z_i$ with $i\leq j-1$, $H^{(j+1)}z_i
=H^{(j)}z_i+\bigl(w_j-H^{(j)}z_j\bigr)z_j^\top z_i
=H^{(j)}z_i.$
Classical inverse-secant algebra \citep{Broyden1965,FangSaad2009} proves its least-change property that \eqref{eq:inner-update} minimizes
the Frobenius norm of the change to $H^{(j)}$ subject to the new secant
condition $H^{(j+1)}z_j=w_j$.

\setcounter{algocf}{\value{algorithm}}

\begin{algorithm2e}[htbp]
\DontPrintSemicolon
\caption{\textsc{Inner-Learning}}
\label{alg:inner}

\KwRequire{
Fixed $J_k$, nonzero $F(v_k)$, $H_{k}$,
and tolerance $0<\eta<1$.
}

Set $j=0$, $H^{(0)}=H_k$ and
$U^{(0)},Q^{(0)},S^{(0)}$ empty.\;

Set $p^{(0)}=-H^{(0)}F(v_k)$.\;

\While{$\norm{F(v_k)+J_kp^{(j)}}_2>\eta\norm {F(v_k)}_2$}{

    \If{$\left(I-U^{(j)}(U^{(j)})^\top \right)p^{(j)}=0$}{
        Choose
        $\ell\in\argmax_{1\le i\le n}\norm{\left(I-U^{(j)}(U^{(j)})^\top \right)e_i}_2$
        and set $u_j = \left(I-U^{(j)}(U^{(j)})^\top \right)e_\ell/\norm{\left(I-U^{(j)}(U^{(j)})^\top \right)e_\ell}_2$.\;
    }
    \Else
    {
    $u_j=\left(I-U^{(j)}(U^{(j)})^\top \right)p^{(j)}/\norm{\left(I-U^{(j)}(U^{(j)})^\top \right)p^{(j)}}_2$. \;
    }

    Compute $g_j=J_ku_j$.\;


    Set $z_j,w_j,H^{(j+1)}$ according to
    \eqref{eq:orthogonalize}--\eqref{eq:inner-update}
    and $U^{(j+1)}=[U^{(j)},u_j]$,
    $Q^{(j+1)}=[Q^{(j)},z_j]$,
    $S^{(j+1)}=[S^{(j)},w_j]$.\;

    $j\leftarrow j+1$ and $p^{(j)}=-H^{(j)}F(v_k)$.\;
}

\KwRet{$(H^{(j)},p^{(j)})$}.\;

\end{algorithm2e}

\setcounter{algorithm}{\value{algocf}}


\section{Global Linear Convergence}\label{sec:global}
In this section, we first prove that the inner learning procedure introduced in Section \ref{sec:completion} satisfies the forcing test in \eqref{eq:forcing}. Then, we bound the certified direction defined by (\ref{eq:correction}) by a constant multiple of the Bellman residual in \eqref{eq:merit}, while further showing that the merit function admits a quadratic upper bound along sufficiently short trial steps, thereby ensuring finite termination of the Armijo backtracking line search \eqref{eq:armijo}. Together, these estimates yield a uniform positive lower bound on the accepted step sizes, based on which Theorem~\ref{thm:global} establishes global convergence
with a linear rate of decay of the Bellman residual in $\ell_2$ norm. 

\begin{lemma}\label{lem:completion}
Consider Algorithm~\ref{alg:inner}.
Starting from any $H^{(0)} \in \mathbb{R}^{n \times n}$, in the $j$-th inner iteration,
(\ref{eq:invariant}), orthonormality of $U^{(j)}$ and $Q^{(j)}$,
$\range(S^{(j)})=\range(U^{(j)})$, and $\bigl\|g_j-Q^{(j)}(Q^{(j)})^\top g_j\bigr\|_2\ge s_0>0$ hold. Moreover, in exact arithmetic, the forcing test (\ref{eq:forcing}) holds after at most $n$ learning queries (equivalently, $n$ inner iterations).
Including forcing tests, the procedure requires at most $2n+1$ products with $J_k$ given by Algorithm \ref{alg:outer}.
\end{lemma}
\begin{proof}
We first prove by induction that every required update in the inner learning procedure is well defined and preserves the following properties: $U^{(j)}$ and $Q^{(j)}$ have orthonormal
columns,
\begin{equation}
 J_kS^{(j)}=Q^{(j)},
 H^{(j)}Q^{(j)}=S^{(j)},~ \mathrm{and} ~
 \range(S^{(j)})=\range(U^{(j)}). \label{eq:range_relation}
\end{equation}
We then further use these properties to establish bounds on the number of iterations required by the inner learning procedure and the number of products with $J_k$.

When $j=0$, the three matrices $U^{(0)},Q^{(0)}$ and $S^{(0)}$ are empty with their column spaces being $\{0\}$, thus the displayed relations in \eqref{eq:range_relation} trivially
hold, and $H^{(0)}$ is finite by assumption. Suppose that these properties
hold after $j<n$ queries and that the proposal
$p^{(j)}=-H^{(j)}F(v_k)$ still fails the forcing test
\eqref{eq:forcing}.
The next three parts show that the next update is
well defined and preserves all the properties in
the $(j+1)$-th iteration.

\emph{Direction selection.}
We first show that a new unit query direction $u_j$
orthogonal to all previous query directions is available.
Let $\Pi_j=I-U^{(j)}(U^{(j)})^\top$.
Since $U^{(j)}$ has $j$ orthonormal columns,
$\Pi_j$ is an orthogonal projector of rank $n-j$.
If $\Pi_jp^{(j)}\ne0$, normalization as specified in
Algorithm~\ref{alg:inner} gives a valid query direction.
Otherwise, the coordinate alternative is available because
\begin{equation}
 \sum_{\ell=1}^n\norm{\Pi_je_\ell}_2^2
 =\operatorname{tr}(\Pi_j^\top\Pi_j)
 =\operatorname{tr}(\Pi_j)=n-j>0.
 \label{eq:coordinate-available}
\end{equation}
Thus either choice gives $\norm{u_j}_2=1$ and
$(U^{(j)})^\top u_j=0$.

\emph{Well-defined normalization.}
The induction hypotheses $J_kS^{(j)}=Q^{(j)}$ and
$\range(S^{(j)})=\range(U^{(j)})$ imply that
\begin{equation}
 \range(Q^{(j)}) = \range(J_kS^{(j)})
 =J_k\range(S^{(j)})
 =J_k\range(U^{(j)}).
 \label{eq:image-span}
\end{equation}
Since $Q^{(j)}$ has orthonormal columns, $\bigl\|g_j-Q^{(j)}(Q^{(j)})^\top g_j\bigr\|_2$
is the distance of $J_ku_j$ to this image space specified in \eqref{eq:orthogonalize}.
Equations~\eqref{eq:matrix-bounds} and
\eqref{eq:image-span} give
\begin{eqnarray}
 \bigl\|g_j-Q^{(j)}(Q^{(j)})^\top g_j\bigr\|_2 &=& \min_{q \in \range (Q^{(j)})} \|g_j - q\|_2
 = \min_{a\in\R^j} \|g_j - J_kU^{(j)}a\|_2\\
 &=& \min_{a\in\R^j}\norm{J_k(u_j-U^{(j)}a)}_2\ge s_0\min_{a\in\R^j}\norm{u_j-U^{(j)}a}_2 \label{eq:Jks} \\
  &=& s_0 \min_{a\in\R^j}\sqrt{\|u_j\|^2_2 - 2a^\top (U^{(j)})^\top u_j + a^\top (U^{(j)})^\top U^{(j)} a} \notag\\
 &=&
 s_0\min_{a\in\R^j}\sqrt{1+\norm a_2^2}  = s_0>0\label{eq:uor}
\end{eqnarray}
where (\ref{eq:Jks}) follows from (\ref{eq:matrix-bounds}), and (\ref{eq:uor}) uses the orthonormality of
$U^{(j)}$, the orthogonality of $u_j$ to its columns,
and $\norm{u_j}_2=1$.
Hence \eqref{eq:orthogonalize} and
\eqref{eq:inverse-image} are well defined.

\emph{Preservation of the induction invariants.}
Orthogonal projection gives $(Q^{(j)})^\top z_j=0$
and $z_j^\top z_j=1$.
From \eqref{eq:inverse-image} and the induction
hypothesis $J_kS^{(j)}=Q^{(j)}$,
\begin{equation}
 J_kw_j
 =\frac{J_ku_j-J_kS^{(j)}(Q^{(j)})^\top g_j}{\bigl\|g_j-Q^{(j)}(Q^{(j)})^\top g_j\bigr\|_2}
 =\frac{g_j-Q^{(j)}(Q^{(j)})^\top g_j}{\bigl\|g_j-Q^{(j)}(Q^{(j)})^\top g_j\bigr\|_2}
 =z_j.
 \label{eq:new-exact-image}
\end{equation}
Using $H^{(j)}Q^{(j)}=S^{(j)}$ and the orthogonality
relations above, the model update
\eqref{eq:inner-update} satisfies
\begin{equation}
 \begin{aligned}
 H^{(j+1)}Q^{(j)}
 &=H^{(j)}Q^{(j)}
   +(w_j-H^{(j)}z_j)z_j^\top Q^{(j)}
 =S^{(j)},\\
 H^{(j+1)}z_j
 &=H^{(j)}z_j
   +(w_j-H^{(j)}z_j)z_j^\top z_j
 =w_j.
 \end{aligned}
 \label{eq:secants-preserved}
\end{equation}
Appending the new columns yields
\begin{equation}
 \begin{aligned}
 J_kS^{(j+1)}
 &=[J_kS^{(j)},J_kw_j]
 =[Q^{(j)},z_j]
 =Q^{(j+1)},\\
 H^{(j+1)}Q^{(j+1)}
 &=[H^{(j+1)}Q^{(j)},H^{(j+1)}z_j]=[S^{(j)},w_j]
 =S^{(j+1)}.
 \end{aligned}
 \label{eq:augmented-invariants}
\end{equation}
Orthonormality of the augmented $U^{(j+1)}$ and
$Q^{(j+1)}$ follows from the construction of $u_j$
and $z_j$.
To check the remaining span relation, using both
directions of \eqref{eq:inverse-image} yields
\begin{equation}
\begin{aligned}
w_j
&=
\frac{u_j-S^{(j)}(Q^{(j)})^\top g_j}
     {\|g_j-Q^{(j)}(Q^{(j)})^\top g_j\|_2}, \quad
u_j
=
\|g_j-Q^{(j)}(Q^{(j)})^\top g_j\|_2\,w_j
+S^{(j)}(Q^{(j)})^\top g_j.
\end{aligned}
\label{eq:primal-span-relations}
\end{equation}
Consequently,
 \begin{align}
 \range(S^{(j+1)})
 &=\range[S^{(j)},w_j]
 =\range[S^{(j)},u_j] \label{eq:primal-span-preserved}\\
 &=\range[U^{(j)},u_j]
 =\range(U^{(j+1)}), \label{eq:induct3}
 \end{align}
where the equality \eqref{eq:induct3} between the spans of $[S^{(j)},u_j]$ and $[U^{(j)},u_j]$ uses the induction hypothesis
$\range(S^{(j)})=\range(U^{(j)})$.
As all induction invariants hold in the $(j+1)$-th iteration, \eqref{eq:range_relation} follows.

\emph{Termination.} Assume that the inner learning procedure has not stopped before $n$ queries,
$Q^{(n)}$ is a square orthogonal matrix.
The two inverse relations constructed by \eqref{eq:invariant} give
\begin{equation}
 H^{(n)}Q^{(n)}
 =S^{(n)}
 =J_k^{-1}Q^{(n)}.
 \label{eq:full-basis-agreement}
\end{equation}
Multiplying their difference by $(Q^{(n)})^\top$ and using $Q^{(n)}(Q^{(n)})^\top = I$ prove
\begin{equation}
 H^{(n)}-J_k^{-1}
 =(H^{(n)}-J_k^{-1})Q^{(n)}(Q^{(n)})^\top
 =0.
 \label{eq:full-inverse}
\end{equation}
The resulting proposal $p^{(n)}=-H^{(n)}F(v_k)$ then must satisfy
\begin{equation}
 \norm{F(v_k)+J_kp^{(n)}}_2
 =\norm{F(v_k)-J_kH^{(n)}F(v_k)}_2
 =\norm{F(v_k)-J_kJ_k^{-1}F(v_k)}_2
 =0,
 \label{eq:exact-policy-residual}
\end{equation}
which uniformly satisfies the forcing test \eqref{eq:forcing}, so no
further learning query is needed.

Finally, if $m\le n$ queries are taken, constructing
$J_ku_j$ costs $m$ products.
Testing the initial proposal and each updated proposal costs at most $m+1$ more. Therefore, the procedure requires at most $2n+1$ products with $J_k$.
\end{proof}

It is worth noting that the inner learning procedure stops as soon as the forcing test \eqref{eq:forcing} is satisfied. As will be shown in Section \ref{sec:experiments}, on public datasets, the inner learning procedure typically performs well with far fewer than $n$ learning queries. 

\begin{lemma}
\label{lem:sign}
Let $v_k \in \mathbb{R}^n$ be a non-optimal value vector in the $k$-th outer iteration, with $b_k = \sgn(F(v_k))$ denoted as the sign vector of the residual and the Bellman residual $F(v_k)$ defined in (\ref{eq:merit}).
Then $b_k$ satisfies
\begin{equation}
 M(v_k; b_k)\ge(1-\gamma)\norm {F(v_k)}_1>0.
 \label{eq:sign-margin}
\end{equation}
Consequently, the direction correction in (\ref{eq:correction}) is well defined.
For any proposal $p_k$ satisfying (\ref{eq:forcing}),
the direction $d_k$ constructed by (\ref{eq:correction}) is nonzero and satisfies
\begin{equation}
 D(v_k;d_k)\le-\theta \|F(v_k)\|_2^2,
 \qquad
 \norm {d_k}_2\le\kappa \|F(v_k)\|_2,
 \label{eq:direction}
\end{equation}
where the uniform bound of the second inequality in (\ref{eq:direction}) is given by $ \kappa
 :=
 \frac{1+\bar\eta}{s_0}
 +\frac{\sqrt n}{1-\gamma}
 \left(
    \frac{L(1+\bar\eta)}{s_0}+\theta
 \right)$ with $s_0$ and $L$ defined in (\ref{eq:matrix-bounds}), and $\bar \eta$ and $\theta$ specified in (\ref{eq:forcing}) and (\ref{eq:correction}), respectively.
\end{lemma}
\begin{proof}
For any row-wise stochastic mixture, $\norm{P_\lambda b}_\infty\le1$.
Since $F(v)^\top b=\norm {F(v)}_1$, this gives
\begin{equation}
 \begin{aligned}
 F(v)^\top J_\lambda b
 =\norm {F(v)}_1-\gamma F(v)^\top P_\lambda b\ge\norm {F(v)}_1-\gamma\norm {F(v)}_1\norm{P_\lambda b}_\infty
 \ge(1-\gamma)\norm {F(v)}_1.
 \end{aligned}
 \label{eq:sign-proof}
\end{equation}
Taking the minimum over the envelope proves \eqref{eq:sign-margin}. Define $r=F(v)+J_\pi p$. By
\eqref{eq:matrix-bounds} and \eqref{eq:forcing},
\begin{equation}
 \norm p_2
 =\norm{J_\pi^{-1}(r-F(v))}_2
 \le s_0^{-1}(\norm r_2+ \|F(v)\|_2)
 \le \frac{1+\bar\eta}{s_0} \|F(v)\|_2.
 \label{eq:proposal-bound}
\end{equation}
If the proposal $p$ satisfies the first row of \eqref{eq:correction}, $d=p$ satisfies $D(v;d)\le-\theta \|F(v)\|_2^2$. Moreover, it follows from \eqref{eq:matrix-bounds} and (\ref{eq:D}) that 
\begin{equation}
    D(v;p) \leq L\|F(v)\|_2\norm p_2 
    \le L\frac{1+\bar\eta}{s_0} \|F(v)\|_2^2 \label{eq:pb}
\end{equation}
where the second inequality in (\ref{eq:pb}) follows from the fact that $s_0\| p\|_2 \leq (1 + \bar\eta)\| F(v) \|_2$ implied by (\ref{eq:matrix-bounds}) and (\ref{eq:forcing}). Therefore, in this case $\| d\|_2 = \|p \|_2 \leq \frac{1+\bar\eta}{s_0} \|F(v)\|_2$, satisfying (\ref{eq:direction}).

If $p$ fails the first row of \eqref{eq:correction}, then $d=p-tb$, where
$t=(D(v;p)+\theta \|F(v)\|_2^2)/M(v;b)>0$. For every admissible mixture,
\begin{equation}
 F(v)^\top J_\lambda d
 =F(v)^\top J_\lambda p-tF(v)^\top J_\lambda b
 \le D(v;p)-tM(v;b)
 =-\theta \|F(v)\|_2^2.
 \label{eq:corrected-pairing}
\end{equation}
Taking the maximum over the envelope proves the descent bound in \eqref{eq:direction}.

For $\| d \|_2$, \eqref{eq:pb}
implies $D(v;p)\le L\|F(v)\|_2\norm p_2\le L\frac{1+\bar\eta}{s_0} \|F(v)\|_2^2$. Hence, by (\ref{eq:sign-proof}) we further bound $0<t
 \le\frac{(L(1+ \bar \eta)+s_0\theta)\|F(v)\|_2^2}{s_0(1-\gamma)\norm {F(v)}_1}
 \le\frac{L(1+ \bar \eta)+s_0\theta}{s_0(1-\gamma)}\|F(v)\|_2$, where the last inequality uses $\norm {F(v)}_1\ge \|F(v)\|_2$. Consequently,
\begin{equation}
 \begin{aligned}
 \norm d_2
 &\le\norm p_2+t\norm b_2
 \le \left [\frac{1+\bar\eta}{s_0}
 +\frac{\sqrt n}{1-\gamma}
 \left(
    \frac{L(1+\bar\eta)}{s_0}+\theta
 \right) \right ] \|F(v)\|_2 = \kappa \|F(v)\|_2,
 \end{aligned}
 \label{eq:correction-length}
\end{equation}
where the last inequality in \eqref{eq:correction-length} uses $\norm b_2\le\sqrt n$.
Finally, $D(v;0)=0$ rules out $d=0$ in the strictly negative descent bound, which completes the proof.
\end{proof}

The following lemma establishes an upper bound on the step size in the directional search defined in \eqref{eq:armijo}, such that all greedy actions remain in the near-active set and the merit function admits a corresponding upper bound.
\begin{lemma}\label{lem:capture}
Consider the non-optimal value vector $v_k$ in the $k$-th outer iteration and suppose the current direction $d_k$ satisfies (\ref{eq:direction}). Let $a_{\rm cap}=\frac{\tau}{2\gamma \kappa \sqrt n}$ with $\tau \in (0, 1-\theta)$, where $\theta$ and $\kappa$ are defined in (\ref{eq:correction}) and (\ref{eq:direction}), respectively. Then,
for any $t \in [0, a_{\rm cap}]$, every action greedy at $v_k + td_k$ belongs to its
near-active set at $v_k$, and 
\begin{equation}
 \Phi(v_k + \alpha d_k)
 \le\Phi(v_k)+\alpha D(v_k; d_k)+\tfrac12L^2\alpha^2\norm {d_k}_2^2, \quad \forall \alpha \in[0, a_{\rm cap}].
 \label{eq:merit-expansion}
\end{equation}
\end{lemma}
\begin{proof}
Let $a$ and $a_0$ be the greedy actions at $v_k+td_k$ and $v_k$ in the state $i$, respectively. It follows from \eqref{eq:q_definitions} that $(Tv_k)_i-q_i^a(v_k) \ge 0$ and
\begin{align}
 (Tv_k)_i-q_i^a(v_k)
 &=q_i^{a_0}(v_k)-q_i^a(v_k) = [q_i^{a_0}(v_k) - q_i^{a_0}(v_k + td_k)] + [q_i^a(v_k + td_k) - q_i^a(v_k)] \notag\\
 &+ [q_i^{a_0}(v_k + td_k) - q_i^a(v_k + td_k)] \notag\\
 & = \gamma t(P_i^a-P_i^{a_0})d_k + q_i^{a_0}(v_k+td_k)-q_i^a(v_k+td_k) \label{eq:greedy_q}\\
 &\le\gamma t\norm{P_i^a-P_i^{a_0}}_1\norm {d_k}_\infty \label{eq:tine} \\
 &\le2\gamma t \kappa \|F(v_k)\|_2, \label{eq:action-gap-capture}
\end{align}
where \eqref{eq:greedy_q} utilizes \eqref{eq:q_definitions}, (\ref{eq:tine}) follows from $q_i^a(v_k+td_k)\ge q_i^{a_0}(v_k+td_k)$, \eqref{eq:action-gap-capture} uses the fact that $P_i^a$ and $P_i^{a_0}$ are both transition probability rows and \eqref{eq:direction}. Thus, if 
$t\le a_{\rm cap}$, the term in \eqref{eq:action-gap-capture} is at most
$\varepsilon(v_k)$ defined in \eqref{eq:width}, proving that any greedy action over $[0, a_{\mathrm{cap}}]$ belongs to $\A_i^\varepsilon(v_k)$ defined in \eqref{eq:width}.

For $t \in [0, a_{\mathrm{cap}}]$, $F(v_k+td_k)$ is continuous and piecewise affine. Its
derivative is $J_{\pi(t)}d_k$ almost everywhere for a greedy policy
$\pi(t)$. Persistent ties have equal derivatives along the segment. Moreover,
for $t \in [0, a_{\mathrm{cap}}]$, the membership  proved above gives
$F(v_k)^\top J_{\pi(t)}d_k\le D(v_k;d_k)$ by the definition of $D(v_k;d_k)$ in (\ref{eq:D}). Consequently,
 \begin{align}
 \frac{d}{dt}\Phi(v_k+td_k)
 &=F(v_k+td_k)^\top J_{\pi(t)}d_k =F(v_k)^\top J_{\pi(t)}d_k
     +(F(v_k+td_k)-F(v_k))^\top J_{\pi(t)}d_k \notag\\
 &\le D(v_k;d_k)
     +\norm{F(v_k+td_k)-F(v_k)}_2\,\norm{J_{\pi(t)}d_k}_2 \le D(v_k;d_k)+L^2t\norm {d_k}_2^2,  \label{eq:along-segment-slope}
 \end{align}
where the second inequality in \eqref{eq:along-segment-slope} comes from \eqref{eq:matrix-bounds} and \eqref{eq:two-bounds}. Absolute continuity now permits:
\begin{align}
 \Phi(v_k+\alpha d_k)-\Phi(v_k)
 =\int_0^\alpha\frac{d}{dt}\Phi(v_k+td_k)\,dt &\le\int_0^\alpha\bigl(D(v_k;d_k)+L^2t\norm {d_k}_2^2\bigr)\,dt \notag\\
 &=\alpha D(v_k;d_k)+\tfrac12L^2\alpha^2\norm {d_k}_2^2, \label{eq:integrated-merit}
\end{align}
which claims \eqref{eq:merit-expansion}.
\end{proof}

\begin{remark}[Feasibility of Armijo backtracking line search]
The merit function $\Phi$ need not be differentiable at policy boundaries, as typically assumed in the standard analysis of the Armijo condition. Nevertheless, Lemma~\ref{lem:capture} provides an explicit range on the step size required by the backtracking line search. For a certified direction at a
non-terminal iterate, \eqref{eq:direction} gives $D(v;d)<0$ and $d\ne0$.
Consequently, (\ref{eq:armijo}) holds for $\alpha \in (0, \bar\alpha(v,d)]$, with $\bar\alpha(v,d)$ given by
\begin{equation}\label{eq:alphabound}
 \bar\alpha(v,d)
 =\min\left\{1,a_{\rm cap},
       \frac{-2(1-\sigma)D(v;d)}{L^2\norm d_2^2}\right\}
 >0,
\end{equation}
since
\eqref{eq:merit-expansion} and the last term in the minimum in (\ref{eq:alphabound}) give
\[
 \begin{aligned}
 \Phi(v+\alpha d)
 \le\Phi(v)+\alpha D(v;d)
               +\tfrac12\alpha^2L^2\norm d_2^2
               &\le\Phi(v)+\alpha D(v;d)
               -(1-\sigma)\alpha D(v;d)\\
               &=\Phi(v)+\sigma\alpha D(v;d),
 \end{aligned}
\]
which is exactly \eqref{eq:armijo}. The restriction $\alpha\le a_{\rm cap}$ is essential, as it keeps every encountered greedy action within the corresponding current near-active set. Thus, the lemma provides an explicit sufficient range for the step size and, in particular, establishes the feasibility of the Armijo  line search. Note that \eqref{eq:alphabound} gives a sufficient condition rather than an upper bound that restricts the step size. In practical implementation, backtracking is therefore retained to test larger step sizes, including the unit step required for the local superlinear convergence analysis. 
\end{remark}
We now present the global linear convergence result of the proposed Algorithm ~\ref{alg:outer}. To this end, define 
\begin{equation}
 a_*:=\min\left\{a_{\rm cap},
              \frac{2(1-\sigma)\theta}{L^2\kappa^2}\right\},
 \label{eq:global-constants}
\end{equation}
where $\theta$, $\sigma$, $\kappa$ and $a_{\rm cap}$ are defined in (\ref{eq:correction}), (\ref{eq:armijo}), \eqref{eq:direction}, and Lemma \ref{lem:capture}, respectively.
\begin{theorem}\label{thm:global}
Consider any initial value vector $v_0 \in \mathbb{R}^n$ and the inverse model $H_0 \in \mathbb{R}^{n \times n}$, with $a_\ast$ and $\beta$ defined in (\ref{eq:global-constants}) and (\ref{eq:armijo}), respectively. 
At each non-terminal iteration $k$, the Armijo backtracking line search uses at most $m^\ast := 1+\left\lceil\frac{\log a_*}{\log\beta}\right\rceil$
trials and accepts $\alpha_k = \beta^{\ell_k-1} \in [\underline\alpha, 1]$ as the step size, where $\ell_k$ denotes the number of trials and $\underline\alpha := \beta a_* >0$. With $\chi=\sqrt{1-2\sigma\theta\underline\alpha}\in(0,1)$, it holds that 
\begin{equation}
\norm{F(v_{k+1})}_2\le\chi\norm{F(v_k)}_2.
 \label{eq:global-rate}
\end{equation}
and 
\begin{equation}
 \norm{v_k-v^*}_\infty
 \le\frac{\chi^k\norm{F(v_0)}_2}{1-\gamma},
 \label{eq:global-value-rate}
\end{equation}
where $v^\ast$ is the optimal value vector. 
\end{theorem}
\begin{proof}
We divide the proof into two parts. The first establishes that the line search terminates within $m^\ast$ trials, while the second establishes the global linear convergence of the value vector.

\emph{Finite termination of the Armijo line search.}
For any \(0<\alpha\le a_*\), Lemma~\ref{lem:capture} applies and
 \begin{align}
 \tfrac12L^2\alpha^2\norm {d_k}_2^2 
\le\tfrac12L^2\kappa^2\alpha^2\|F(v_k)\|_2^2 \le(1-\sigma)\alpha\theta \|F(v_k)\|_2^2 \le-(1-\sigma)\alpha D(v_k;d_k),  \label{eq:usdi}
\end{align}
where the first inequality in \eqref{eq:usdi} uses (\ref{eq:direction}), the second inequality comes from the second term in \eqref{eq:global-constants} and the third comes from \eqref{eq:direction}. Substituting \eqref{eq:usdi} into \eqref{eq:merit-expansion} yields the following Armijo-type sufficient decrease condition required by \eqref{eq:armijo}:
\begin{equation}
 \Phi(v_k+\alpha d_k)
 \le\Phi(v_k)+\alpha D(v_k;d_k)-(1-\sigma)\alpha D(v_k;d_k)
 =\Phi(v_k)+\sigma\alpha D(v_k;d_k).
 \label{eq:armijo-guarantee}
\end{equation}
Thus \eqref{eq:armijo} holds for $\alpha \in (0, a_*]$ . Since the Armijo backtracking line search starts from a step size of $1$ and reduces it by a factor of $\beta$ at each trial, the line search terminates after at most $1+\left\lceil\frac{\log(a_*)}{\log(\beta)}\right\rceil$
trials. In particular, let $\ell_k$ denote the number of trials required by the line search at the $k$-th outer iteration. Then, the accepted step size is given by $\alpha_k=\beta^{\ell_k-1}\in[\beta a_*,1].$


\emph{Global linear convergence.}
The accepted step size, together with the fact that $d_k\neq 0$ by Lemma~\ref{lem:sign}, implies that $s_k=\alpha_k d_k\neq 0$. By
\eqref{eq:two-bounds},
\begin{equation}
 \norm{y_k}_\infty
=\norm{F(v_{k+1})-F(v_k)}_\infty\ge(1-\gamma)\norm{s_k}_\infty>0.
 \label{eq:trajectory-nonzero}
\end{equation}
Hence \eqref{eq:trajectory} has a positive denominator and produces a well defined $H_{k+1}$. 

The accepted step size and \eqref{eq:armijo} imply that
\begin{eqnarray}
 \Phi(v_{k+1})
 &\le&\Phi(v_k)+\sigma\alpha_kD(v_k;d_k)\\
 &\le&\Phi(v_k)-\sigma\alpha_k\theta\norm{F(v_k)}_2^2\label{eq:merit-contraction}\\
 &\le&(1-2\sigma\theta\underline\alpha)\Phi(v_k), \label{eq:root}
\end{eqnarray}
where (\ref{eq:merit-contraction}) follows from \eqref{eq:direction}.
By the prescribed parameter ranges, we have $0<2\sigma\theta\underline\alpha<1$.
Taking square roots of (\ref{eq:root}) yields \eqref{eq:global-rate}. It further follows from \eqref{eq:error} that
\begin{equation}
 \norm{v_k-v^*}_\infty
 \le\frac{\norm{F(v_k)}_\infty}{1-\gamma}
 \le\frac{\norm{F(v_k)}_2}{1-\gamma}
 \le\frac{\chi^k\norm{F(v_0)}_2}{1-\gamma},
 \label{eq:value-error-chain}
\end{equation}
which proves \eqref{eq:global-value-rate}.
\end{proof}
Theorem \ref{thm:global} proves that the Bellman residual $\|F(v_k)\|_2$ converges Q-linearly to zero, while the value iterates $v_k$ converge R-linearly to the optimal value vector $v^*$. Note that these results do not require the inverse model $H_k$ to be nonsingular or uniformly bounded in norm, as required in \citep{BurdakovKamandi2018,ThemelisPatrinos2019,LiuEtAl2023BlockBroyden}. The next section subsequently studies the local behavior of the proposed algorithm, establishing a superlinear rate.

\section{Local Superlinear Convergence}\label{sec:local}
While the global linear convergence of $v_k\to v^*$ has been established in Section~\ref{sec:global} based on the line search method and the inner learning procedure, quasi-Newton methods are known to achieve local superlinear convergence for nonlinear equations and fixed-point problems \citep{BroydenDennisMore1973,SunHan1997,ThemelisPatrinos2019}. Since the inner learning procedure plays a role analogous to that of a quasi-Newton method in approximating the inverse Jacobian, we show in this section that the proposed algorithm inherits this desirable local convergence behavior.

We first show that the near-active set eventually contains only optimal actions, without requiring the existence of a unique optimal policy. This allows continuous switching among different policy boundaries while still guaranteeing a decrease in the merit function, thereby addressing a key difficulty in decision-making problems involving policy changes.
\begin{lemma}\label{lem:identification}
Consider Algorithm~\ref{alg:outer} with the outer iteration indexed by $k$ and define 
\begin{equation}\label{eq:Astar}
\A_i^*=\argmax_{a\in\A}q_i^a(v^*).
\end{equation}
Then, there is an index $K_0 >0$ such that
\begin{equation}
 \A_i^\varepsilon(v_k)\subseteq\A_i^*
 \quad\text{for every }i \in \mathcal{S}\text{ and }k\ge K_0,
 \label{eq:optimal-containment}
\end{equation}
where $\A_i^\varepsilon(v_k)$ is defined in (\ref{eq:width}).
\end{lemma}
\begin{proof}
Define $g_i^a(v)=(Tv)_i-q_i^a(v)$ as the gap between action $a$ and the greedy action at value vector $v$. If $\mathcal A_i^*=\mathcal A$ for every state $i$, the claim is
immediate. Otherwise, since the action set is finite,
there exists a positive gap $\Delta_*$ separating the optimal and non-optimal state-action pairs, i.e.,
\begin{equation}
 \Delta_*:=\min_{\{(i,a):a\notin\A_i^*\}}g_i^a(v^*)>0.
 \label{eq:strict-gap}
\end{equation}
By continuity of $g_i^a(\cdot)$ and $\varepsilon(\cdot)$, the convergence $v_k\to v^\ast$ established in Theorem~\ref{thm:global}, and the finiteness of the state-action space, there exists some $K_0>0$ such that for all $k \geq K_0$, $\varepsilon(v_k)$ defined in \eqref{eq:width} satisfies $\varepsilon(v_k)<\Delta_*/2$ and $g_i^a(v_k) \ge \Delta_*/2$ for each strictly suboptimal pair $(i,a)$. Consequently, for any such pair, it follows from \eqref{eq:width} that
\begin{equation}
 g_i^a(v_k)\ge\Delta_*/2>\varepsilon(v_k)
 \quad\Longrightarrow\quad
 a\notin\A_i^\varepsilon(v_k),
 \label{eq:gap-exclusion}
\end{equation}
\end{proof}

\begin{lemma}\label{lem:disagreement}
Suppose $\A_i^\varepsilon(v_k)\subseteq\A_i^*$ for every state $i \in \mathcal{S}$ in the $k$-th outer iteration, with $\A_i^*$ defined in (\ref{eq:Astar}). Then, every stationary randomized policy $\lambda$ supported on the near-active sets
satisfies
\begin{equation}
 F(v_k)^\top J_\lambda e_k\ge(1-\tau)\|F(v_k)\|_2^2,
 \label{eq:disagreement}
\end{equation}
where $e_k = v_k -v^*$, $\tau$ and $F(v_k)$ are defined in (\ref{eq:width}) and (\ref{eq:merit}), respectively. 
\end{lemma}
\begin{proof}
For any value vector $v_k$ at outer iteration $k$, with $\pi_k$ the greedy policy, from (\ref{eq:q_definitions}) we have 
\begin{eqnarray}
    F(v_k) =  v_k - Tv_k &=& (I - \gamma P_{\pi_k})v_k - r_{\pi_k} \label{eq:nodet} \\
    & = & (I - \gamma P_{\pi_k})v_k - (I - \gamma P_{\pi_k})v^\ast \label{eq:J_pi_e} \\
    & = & J_{k} e_k \label{eq:combine}, 
\end{eqnarray}
where in (\ref{eq:nodet}) $r_{\pi_k} = [r_1^{a_1}, \cdots, r_n^{a_n}]^\top$ denotes the reward vector, with $a_i\in\argmax_{a\in\A}q_i^a(v_k)$ the chosen greedy action under $\pi_k$ for $i\in\{1,\ldots,n\}$, and \eqref{eq:J_pi_e} follows from the fact that $\pi_k$ is also a greedy policy with respect to the optimal value vector $v^\ast$, under the assumption that $\A_i^\varepsilon(v_k)\subseteq\A_i^*$. For each $a \in \A_i^\varepsilon(v_k)$,  $\A_i^\varepsilon(v_k)\subseteq\A_i^*$ also implies that
$v_i^*=r_i^a+\gamma P_i^av^*$, subtracting $F_i(v_k)$ from both sides gives
 \begin{align}
 e_{k,i}-\gamma P_i^ae_k-F_i(v_k)
 =v_{k,i}-\gamma P_i^av_k-r_i^a-F_i(v_k)=(Tv_k)_i-q_i^a(v_k) := g_i^a(v_k). \label{eq:gap-identity}
 \end{align}
For any stationary randomized policy $\lambda$ supported on the near-active sets at $v_k$, define
$(\delta_{\lambda,k})_i=\sum_a\lambda_{i}^ag_i^a(v_k)$.
\eqref{eq:width} implies that $0\le(\delta_{\lambda,k})_i\le\varepsilon(v_k)$, and it follows from (\ref{eq:gap-identity}) that
\begin{eqnarray}
    (\delta_{\lambda,k})_i
&=& \sum_{a\in\mathcal A_i^\varepsilon(v_k)}
  \lambda_i^a g_i^a(v_k) = \sum_{a\in\mathcal A_i^\varepsilon(v_k)}
  \lambda_i^a\bigl(e_{k,i}-\gamma P_i^ae_k-F_i(v_k)\bigr) \label{eq:form} \\
& = & e_{k,i}-\gamma(P_\lambda)_i e_k-F_i(v_k)=(J_\lambda e_k)_i-F_i(v_k). \label{eq:compact}
\end{eqnarray}
where the first equality in (\ref{eq:compact}) follows from the fact that \(\sum_a\lambda_i^a=1\) and \((P_\lambda)_i=\sum_a\lambda_i^aP_i^a\). Therefore, 
\eqref{eq:compact} gives $J_\lambda e_k=F(v_k)+\delta_{\lambda,k}$, and it follows that
 \begin{align}
 F(v_k)^\top J_\lambda e_k
 &=\|F(v_k)\|_2^2+F(v_k)^\top\delta_{\lambda,k} \ge \|F(v_k)\|_2^2-\norm {F(v_k)}_1\norm{\delta_{\lambda,k}}_\infty \notag\\
 &\ge \|F(v_k)\|_2^2-\sqrt n \|F(v_k)\|_2\,\varepsilon(v_k)
 =(1-\tau)\|F(v_k)\|_2^2,  \label{eq:disagreement-proof}
\end{align}
where \eqref{eq:disagreement-proof} comes from \eqref{eq:width} that $\varepsilon(v_k)=\tau \|F(v_k)\|_2/\sqrt n$.
\end{proof}
The following theorem establishes the local superlinear convergence of the proposed algorithm.
\begin{theorem}\label{thm:local}
Consider Algorithm \ref{alg:outer},  with the outer iteration indexed by $k$, and let $K_0$ be defined in Lemma \ref{lem:identification}. If it does not terminate finitely, then there exists $K_1 > K_0$ such that, for $k \geq K_1$, the direction of line search $d_k$, the proposal $p_k$ and the step size $\alpha_k$ obey $d_k=p_k$ and $\alpha_k=1$. Moreover, it holds that
\begin{equation}
 \norm{v_{k+1}-v^*}_2
 \le\frac{cL^{1+\nu}}{s_0}\norm{v_k-v^*}_2^{1+\nu},
 \label{eq:local-rate}
\end{equation}
where $L$ and $s_0$ are defined in (\ref{eq:matrix-bounds}), while $c$ and $\nu$ are defined in (\ref{eq:forcing}).
\end{theorem}
\begin{proof}

Define $r_k=F(v_k)+J_{k} p_k$. As  $\A_i^\varepsilon(v_k)\subseteq\A_i^*$
for every $i \in \mathcal{S}$ and $k\ge K_0$ by Lemma \ref{lem:identification}, it follows from \eqref{eq:combine} in Lemma \ref{lem:disagreement} that
\begin{equation}
 r_k=J_{k}(e_k+p_k),\qquad p_k=-e_k+J_{k}^{-1}r_k, \qquad \forall k \geq K_0.
 \label{eq:local-proposal}
\end{equation}
From \eqref{eq:matrix-bounds} and \eqref{eq:forcing},
\begin{equation}
 \norm{e_k+p_k}_2
 =\norm{J_{k}^{-1}r_k}_2
 \le s_0^{-1}\norm {r_k}_2
 \le s_0^{-1}\eta_k \|F(v_k)\|_2.
 \label{eq:error-step}
\end{equation}

For any stationary randomized policy supported on the near-active set, we have
\begin{align}
F(v_k)^\top J_\lambda p_k
&=-F(v_k)^\top J_\lambda e_k
  +F(v_k)^\top J_\lambda J_k^{-1}r_k
  \label{eq:correction-substitution}\\
&\le -(1-\tau)\|F(v_k)\|_2^2
  +\|F(v_k)\|_2\|J_\lambda\|_2
   \|J_k^{-1}\|_2\|r_k\|_2
  \label{eq:correction-estimate}\\
&\le -\left(1-\tau-\frac{L}{s_0}\eta_k\right)
  \|F(v_k)\|_2^2,
  \label{eq:correction-bound}
\end{align}
where \eqref{eq:correction-substitution} follows from \eqref{eq:local-proposal},
\eqref{eq:correction-estimate} follows from \eqref{eq:disagreement} and the
Cauchy--Schwarz inequality together with the induced matrix norm, and \eqref{eq:correction-bound} uses
$\|J_\lambda\|_2\le L$, $\|J_k^{-1}\|_2\le s_0^{-1}$
from \eqref{eq:matrix-bounds}, and $\|r_k\|_2\le\eta_k\|F(v_k)\|_2$ from \eqref{eq:forcing}.
Maximizing over $\lambda$ then yields the same bound for $D(v_k;p_k)$. Since $1-\tau-\theta>0$ by the choice $\tau\in(0,1-\theta)$, it follows that
\begin{equation}\label{eq:evpass}
    D(v_k;p_k)\le -\theta\|F(v_k)\|_2^2
\end{equation}
whenever \( \eta_k\le \frac{s_0}{L}(1-\tau-\theta).\)
This condition eventually holds for all $k \geq \bar{K}$ with $\bar{K} \geq K_0$ by the forcing test \eqref{eq:forcing} and the fact that $\|F(v_k)\|_2\to0$. Therefore, $d_k=p_k$ by \eqref{eq:correction}.

The current greedy policy, viewed as a stationary randomized policy assigning probability one to its chosen action in each state, is supported on the near-active set, i.e., $J_{k}^\top F(v_k) \in \G_\varepsilon(v_k)$. Then it follows from (\ref{eq:D}) that 
\begin{equation}
 D(v_k;p_k)\ge F(v_k)^\top J_{k} p_k
 =F(v_k)^\top(r_k-F(v_k))
 \ge-(1+\eta_k)\|F(v_k)\|_2^2,
 \label{eq:local-support-lower}
\end{equation}
where the final inequality in \eqref{eq:local-support-lower} follows from $\norm {r_k}_2\le\eta_k \|F(v_k)\|_2$.
This lower bound  controls how
small the right-hand side of Armijo condition (\ref{eq:armijo}) can be:
\begin{equation}
 \Phi(v_k)+\sigma D(v_k;p_k)
 \ge\left(\tfrac12-\sigma(1+\eta_k)\right)\|F(v_k)\|_2^2.
 \label{eq:full-armijo-rhs}
\end{equation}
On the other hand, $F(v^*)=0$, \eqref{eq:two-bounds}, and
\eqref{eq:error-step} together give the merit function taking full step obeying 
\begin{equation}
 \Phi(v_k+p_k)
=\tfrac12\norm{F(v_k+p_k)-F(v^*)}_2^2\le\tfrac12L^2\norm{e_k+p_k}_2^2
 \le\tfrac12(L/s_0)^2\eta_k^2\|F(v_k)\|_2^2.
 \label{eq:full-armijo-lhs}
\end{equation}
To make the comparison explicit, let $a_0=1/2-\sigma>0$.
For $\eta_k\le\min\left\{\frac{a_0}{2\sigma},\,\frac{s_0}{L}\sqrt{a_0}\right\},$
the coefficient $\tfrac12-\sigma(1+\eta_k)$ in \eqref{eq:full-armijo-rhs} is at least $a_0/2$,
and the coefficient $\tfrac12(L/s_0)^2\eta_k^2$ in \eqref{eq:full-armijo-lhs} is at most $a_0/2$.
Hence, \eqref{eq:armijo} holds for the full step. Since $\eta_k\to0$, this condition holds for all sufficiently large $k$; in particular, there exists some $\tilde K\geq \bar{K}$ such that the full step is accepted for all $k\geq\tilde K$. Then it holds that $e_{k+1}=e_k+p_k$ for $k \geq \tilde K$. Moreover, since $\|F(v_k)\|_2 \to 0$, (\ref{eq:forcing}) implies that there exists some $\hat K$ such that $\eta_k = c\|F(v_k)\|_2^\nu$ for all $k \geq \hat K$. Combining \eqref{eq:error-step} with the Lipschitz bound $\|F(v_k)\|_2\leq L\|e_k\|_2$ yields, for all $k\geq K_1:=\max\{\tilde K,\hat K\}$,
 \begin{align}
 \norm{e_{k+1}}_2
 \le s_0^{-1}\eta_k \|F(v_k)\|_2
 =\frac{c}{s_0}\|F(v_k)\|_2^{1+\nu}\le\frac{cL^{1+\nu}}{s_0}\norm{e_k}_2^{1+\nu},  \label{eq:rate-chain}
 \end{align}
which proves \eqref{eq:local-rate}.
\end{proof}

The local superlinear convergence follows from the classical inexact-Newton principle: the relative residual of the approximate Newton equation tends to zero as the iterate approaches the optimal value \citep{DemboEisenstatSteihaug1982}. Since the inner learning procedure satisfies the forcing test at every iteration and the line search eventually accepts the full step, local superlinear convergence is guaranteed. Moreover, the result concerns this complete
algorithm, assuming neither a Dennis--Mor\'e
condition nor policy constancy. In the next section, we further discuss the inverse learning errors and the computational complexity brought by LSQN.

\section{Inverse Learning and Computational Complexity}\label{sec:learning}
As shown in Sections \ref{sec:global}--\ref{sec:local}, the algorithm guarantees the globally linear convergence and the locally superlinear convergence results without requiring uniform accuracy of the entire inverse model. Through this section, we analyze the estimation error bound of the inverse model constructed through the inner learning procedure and the inverse Broyden update in (\ref{eq:trajectory}). We further show that, under policy consistency, the inverse Broyden update reduces the inverse-model error. Along this line, we further investigate a policy-consistent variant of Algorithm \ref{alg:outer} and establish a conditional information bound on the cumulative reduction in the inverse-model error. Finally, we establish a worst-case computational complexity bound for achieving a prescribed residual accuracy.

\subsection{Inverse-Model Error Identities}
Fix an invertible policy Jacobian $J$ and define $E=H-J^{-1}$.
For a unit image $z$ with exact inverse image $w=J^{-1}z$, the rank-one update
\eqref{eq:inner-update} gives the next inverse model $H^+$ and yields
\begin{equation}
 H^+-J^{-1}
 =H+(w-Hz)z^\top -J^{-1} = E+(J^{-1}z-Hz)z^\top=E-Ezz^\top=E(I-zz^\top).
 \label{eq:error-projection}
\end{equation}
The two matrices $E(I-zz^\top)$ and $Ezz^\top$ are orthogonal in the
Frobenius inner product. Indeed,
\begin{equation}
 \left\langle E(I-zz^\top),Ezz^\top\right\rangle_F
 =z^\top(I-zz^\top)E^\top Ez=0.
 \label{eq:error-orthogonality}
\end{equation}
Since $\norm{Ezz^\top}_F=\norm{Ez}_2$, the Pythagorean identity yields
 \begin{align}
 \norm{H^+-J^{-1}}_F^2
 &=\norm E_F^2-\norm{Ezz^\top}_F^2 \notag \\
 &=\norm{H-J^{-1}}_F^2-\norm{(H-J^{-1})z}_2^2, \label{eq:learning-energy}
\end{align}
where \eqref{eq:learning-energy} comes from the classical orthogonal least-change identity
\citep{DennisSchnabel1979,GowerRichtarik2017}, showing exactly the reduction in the inverse-model error under the current policy achieved by each learning query.

For an observed pair $(s,y)$ with $y\ne0$, set
$\zeta=y/\norm y_2$ and $\xi=(s-J^{-1}y)/\norm y_2$.
Subtract $J^{-1}$ from
\eqref{eq:trajectory}, with $H$ denoting the incoming model $\widehat H_k$ in \eqref{eq:trajectory}:
\begin{equation}
 H^+-J^{-1}=E+\frac{(s-Hy)y^\top}{y^\top y}=E+\left(\frac{s-J^{-1}y}{\norm y_2}
          -E\zeta\right)\zeta^\top=E(I-\zeta\zeta^\top)+\xi\zeta^\top.
          \label{eq:inconsistent-error}
\end{equation}
The cross term again vanishes since \( \left\langle
 E(I-\zeta\zeta^\top),\xi\zeta^\top
 \right\rangle_F=\zeta^\top(I-\zeta\zeta^\top)E^\top\xi=0.\)
 Using \eqref{eq:learning-energy} for the projected part
therefore gives
 \begin{align}
 \norm{H^+-J^{-1}}_F^2-\norm{H-J^{-1}}_F^2
 &=-\norm{(H-J^{-1})\zeta}_2^2+\norm\xi_2^2\notag \\
 &=-\norm{(H-J^{-1})\zeta}_2^2
   +\frac{\norm{s-J^{-1}y}_2^2}{\norm y_2^2}. \label{eq:trajectory-energy}
\end{align}
Equation \eqref{eq:trajectory-energy} separates the error removed by the update from the error introduced by secant inconsistency
\citep{DennisWalker1985}. If the same policy $\pi_v$ is greedy at both $v$ and $v+s$, then
\begin{equation}
 y=F(v+s)-F(v)
 =[J_{\pi_v}(v+s)-r_{\pi_v}]-[J_{\pi_v} v-r_{\pi_v}]=J_{\pi_v} s.
 \label{eq:compatible-trajectory}
\end{equation}
Thus, the second term in \eqref{eq:trajectory-energy} vanishes when $J = J_{\pi_v}$. 
\eqref{eq:trajectory-energy} explains the role of the retained inverse model in LSQN. It provides a useful initial candidate step, particularly when the policy remains unchanged, while the forcing test determines whether additional learning is needed currently. 

\subsection{A Conditional Information Bound}
In this subsection, we consider a variant that always retains inner learning but applies the inverse Broyden update specified in \eqref{eq:trajectory} only if $\pi_{k+1}=\pi_k$. This ensures that every applied update satisfies an exact secant relation for the current policy, allowing the error reductions to be accumulated. As the policy changes,
the variant simply retains the model obtained from inner learning.

\begin{proposition}\label{prop:variation}
At outer iteration $k$, let $q_k$ sum the nonnegative inverse-error
reductions in (\ref{eq:learning-energy}) from all learning queries and any
permitted inverse Broyden update, each measured against $J_k^{-1}$, the inverse Jacobian in the $k$-th outer iteration.
Then, over any completed iterations $0, \ldots, K$,
\begin{equation}
 \sum_{k=0}^{K}q_k
 \le\left(A_0+\sum_{k=0}^{K-1}\Delta_k\right)^2,
 \label{eq:variation-budget}
\end{equation}
where $A_0=\norm{H_0-J_0^{-1}}_F$ and
$\Delta_k=\norm{J_{k+1}^{-1}-J_k^{-1}}_F$. Moreover, the conclusions of Theorems~\ref{thm:global}
and~\ref{thm:local} hold as well for this variant.
\end{proposition}
\begin{proof}
Define $A_k=\norm{H_k-J_k^{-1}}_F$ and
$B_k=\norm{H_{k+1}-J_k^{-1}}_F$. Thus both errors at iteration $k$
use the same reference inverse, even if the next selected policy changes.
Every applied update has an exact secant for $J_k$ by
\eqref{eq:new-exact-image} or \eqref{eq:compatible-trajectory}. Summing
\eqref{eq:learning-energy} within that iteration gives
\begin{equation}
 B_k^2=A_k^2-q_k,\qquad B_k\le A_k.
 \label{eq:within-policy-energy}
\end{equation}
Changing the reference at the next iterate gives
 \begin{align}
 A_{k+1}=\norm{H_{k+1}-J_{k+1}^{-1}}_F \le\norm{H_{k+1}-J_k^{-1}}_F
      +\norm{J_k^{-1}-J_{k+1}^{-1}}_F
 =B_k+\Delta_k. \label{eq:reference-change}
\end{align}
Induction using $B_k\le A_k$ yields
\begin{equation}
 B_k\le A_k\le A_0+\sum_{j=0}^{k-1}\Delta_j.
 \label{eq:error-budget}
\end{equation}
Now telescope \eqref{eq:within-policy-energy}:
 \begin{align}
 \sum_{k=0}^{K}q_k
 &=A_0^2-B_K^2+\sum_{k=0}^{K-1}(A_{k+1}^2-B_k^2) \le A_0^2+\sum_{k=0}^{K-1}
\bigl((B_k+\Delta_k)^2-B_k^2\bigr) \label{eq:tele_substitue}\\
 &=A_0^2+\sum_{k=0}^{K-1}(2B_k\Delta_k+\Delta_k^2) \le A_0^2+\sum_{k=0}^{K-1}
\left[2\left(A_0+\sum_{j<k}\Delta_j\right)\Delta_k+\Delta_k^2\right] \label{eq:tele_B_k_sub}\\
 &=\left(A_0+\sum_{k=0}^{K-1}\Delta_k\right)^2. \label{eq:telescoping-information}
\end{align}
where the inequality in \eqref{eq:tele_substitue} comes from \eqref{eq:reference-change} and drops $-B_K^2$, \eqref{eq:tele_B_k_sub} uses \eqref{eq:error-budget}, and  \eqref{eq:telescoping-information} utilizes the expansion of the square. When the outer inverse Broyden update is skipped, $H_{k+1}=\widehat H_k$ remains finite. From previous analyses, the inner learning procedure still attains the prescribed forcing
accuracy in finitely many queries, while the direction construction
and Armijo line search remain unchanged. The proofs of
Theorems~\ref{thm:global} and~\ref{thm:local} therefore can also apply to this variant.
\end{proof}

When the policy Jacobian $J$ remains unchanged, the cumulative
reduction in squared inverse-model error cannot exceed the initial
squared error. After a policy change, the retained model may be less accurate for the new policy Jacobian. The bound includes this effect through
$\Delta_k=\norm{J_{k+1}^{-1}-J_k^{-1}}_F$. Overall, the purpose of retaining the inverse model is to use earlier
learning when computing later search directions. For the variant
considered here, the bound relates the accumulated reduction in
inverse-model error to the initial error and changes in $J_k^{-1}$.
It thus helps explain both the value of preserving previous
learning and the additional corrections that policy changes
may require.

\subsection{Worst-Case Computational Cost Per Outer Iteration} \label{subsec:wccc}
We now analyze the computational cost of a single non-terminal outer iteration $k$. Consider an MDP with states number of $n = |\mathcal{S}|$ and $A = |\mathcal{A}|$ the number of actions, $m_k\le n$ the number of learning queries performed by the inner learning, and $\ell_k$ the number of trial steps defined in Theorem \ref{thm:global}. We first consider dense transition matrices and an explicitly stored dense inverse model.

\textit{Policy selection and direction certification.} Computing the action values at \(v_k\), and hence the Bellman residual, the greedy policy, and the near-active sets, requires \(O(An^2)\) operations in total. These quantities can be reused if they are already available from the preceding accepted trial. Constructing the current policy matrix costs at most \(O(n^2)\). Evaluating the near-active descent certificate and constructing the corrected direction require a constant number of action-wise matrix–vector computations, resulting in a cost of \(O(An^2)\). Therefore, the overall computational complexity of policy selection and direction certification is \(T^{\mathrm{PD}}_k := O(An^2)\).

\textit{Inner learning.} The inner learning procedure accounts for the potential cubic cost within an outer iteration. With $m_k$ learning queries, it requires at most $2m_k+1$ policy-Jacobian products from Lemma \ref{lem:completion}, including those used to check the forcing test \eqref{eq:forcing}. Each dense policy-Jacobian product costs \(O(n^2)\). Applying the dense inverse model and performing each rank-one inverse update \eqref{eq:inner-update} also cost \(O(n^2)\). At inner index \(j\), projection and orthogonalization against the previously acquired basis columns cost \(O(nj)\). Summing these operations over learning queries gives
\[
\sum_{j=0}^{m_k-1}O(nj)=O(nm_k^2).
\]
Consequently, the computational complexity of inner learning is
\[
T_k^{\mathrm{inner}}
=
O\!\left((m_k+1)n^2+nm_k^2\right).
\]
Since \(m_k\le n\), its worst-case cost is \(O(n^3)\). This worst case corresponds to learning a full set of \(n\) independent inverse-secant directions. The inner learning procedure terminates earlier whenever the prescribed forcing accuracy is attained. 

\textit{Armijo line search.} Each line search trial requires evaluating the Bellman residual at \(v_k+\alpha d_k\), at a cost of \(O(An^2)\). Thus, the computational cost of the line search is \(T_k^{\mathrm{LS}}=O(\ell_k An^2)\), where \(\ell_k\) denotes the number of trials, as defined in Theorem \ref{thm:global}.

Moreover, choosing \(0<\beta<1\) such that \(\beta\le a_*\) guarantees that the line search terminates within at most two trials, yielding \(T_k^{\mathrm{LS}}=O(An^2)\). However, this choice of \(\beta\) does not necessarily improve overall computational efficiency, as a small \(\beta\) may lead to shorter accepted steps and consequently increase the number of outer iterations. 

\textit{Per-iteration Cost}. After the Armijo line search accepts a step size $\alpha_k$, the algorithm updates $v_{k+1} = v_k + \alpha_kd_k$ and performs the rank-one inverse update, costing $O(n^2)$. Combining the previous items under the choice $\beta\le a_*$ yields the overall cost as
\[T_k=O\!\left(An^2
+(m_k+1)n^2+nm_k^2\right).
\]
Using the fact that \(m_k\le n\), we obtain the worst case computational cost:
\begin{equation}
T_k= O(n^3 + An^2), \label{eq:overall_cost}
\end{equation}
Note that the worst-case bound does not mean that every outer iteration incurs this computational cost. For instance, if the actual number of learning queries satisfies \(m_k=O(n^\rho)\) for some \(0\leq\rho<1\), then the inner-learning cost becomes
\[
T_k^{\mathrm{inner}}
=O\!\left(n^{2+\rho}+n^{1+2\rho}\right)
=O(n^{2+\rho}).
\]
In Section \ref{sec:experiments}, we illustrate that under all of our settings and examples, the realized learning query rarely reaches $n$ and in most cases usually behaves better than other baselines in runtime.

\section{Numerical Experiments}\label{sec:experiments}
In this section, we evaluate LSQN on four MDP families and public benchmark environments. We first describe the test
problems and solver settings, and then compare LSQN with nine baselines in terms of outer iteration counts and runtime.
We further examine the learning-query counts on the public
benchmarks to assess how much inner learning is needed in
practice relative to the worst-case bound of $n$ queries.
Finally, a separate synthetic FrozenLake experiment investigates
how the mean learning-query count grows with the number of
states, providing an empirical estimate of the associated
computational cost. All numerical results reported in this paper were obtained on an Apple M2 CPU running macOS. We also conducted the experiments on an Intel Core i9-14900K CPU running Linux and obtained similar results. The corresponding results and details of both computing environments are provided in the same repository. The code is available at \url{https://github.com/wang-zhizhou/LSQN.git}.

\subsection{Test Problems and Experimental Setup} \label{subsec:experiment_setup}
For a general test, we use four model-based MDP families:
a controlled queue, a nearly decomposable clustered chain, a sparse Garnet
model, and a dense latent-mixture model, where each is employed with four actions, and Appendix~\ref{app:models} specifies their construction. 
Crossing these families, we experiment with $n \in \{512,1024,2048,4096\}$, using $\gamma \in \{0.99,0.999\}$, and giving a total of 32 instances. 

For comparison, we evaluate 9 baseline methods alongside LSQN, and three runs for each instance are considered to reduce the timing variability. In specific, the baselines are VI, exact PI, MPI with four or twenty evaluation sweeps, three QPI variants, inexact GMRES policy iteration (iGMRES-PI), and stabilized type-I Anderson acceleration (AA-I-S-$m$). For three types of QPI, QPI-U uses the uniform prior and original residual safeguard, while QPI-BT implements (16)--(17) of \citet{KolarijaniEsfahani2023QPI}. Additionally, QPI-S adds their policy-dependent secant constraints in equation (18), using a rank-revealing constrained update and a Woodbury inverse. For the type-I Anderson acceleration on VI \citep{ZhangODonoghueBoyd2020}, we just consider its stabilized variant with VI fallback for comparison. 
For inexact PI, \citet{Gargiani2024IPI} apply GMRES \citep{Gargiani2023GMRES} to solve the policy-evaluation linear system approximately.

In addition, in subsection \ref{subsec:public-inner-learning}, we illustrate that under LSQN, the learning queries in the inner learning procedure is usually far below than $n$, thus reducing the computational work relative to the $O(n^3)$ cost of full inner learning. We use the public benchmark environments FrozenLake, Taxi, and CliffWalking provided by \cite{Towers2024Gymnasium} to test our algorithm.

All solvers start from the same constant reward lower bound and target
$\norm{F(v_k)}_\infty\le10^{-8}$. In our algorithm, we set the initial inverse model as the identity matrix $H_0 = I$, where common parameters are chosen as $\tau = 0.1$, $\theta = 0.5$, $\beta = 0.5$, the forcing test parameters $\bar\eta=c=\nu=0.5$, and the Armijo coefficient $\sigma = 10^{-4}$. Through the test, timings include initialization, policy-matrix assembly, linear
algebra, and rejected trials, while excluding problem generation and output writing. Moreover, a method is deemed successful on a given instance only if all three runs are successful.

\subsection{Runtime and Residual Reduction} \label{subsec:runtime}
In this subsection, we compare our algorithm with other baselines, from the perspective of both the number of outer iterations and runtime. 

\begin{figure}[t]
    \centering

    \begin{minipage}[b]{0.48\linewidth}
        \centering
        \includegraphics[width=\linewidth]{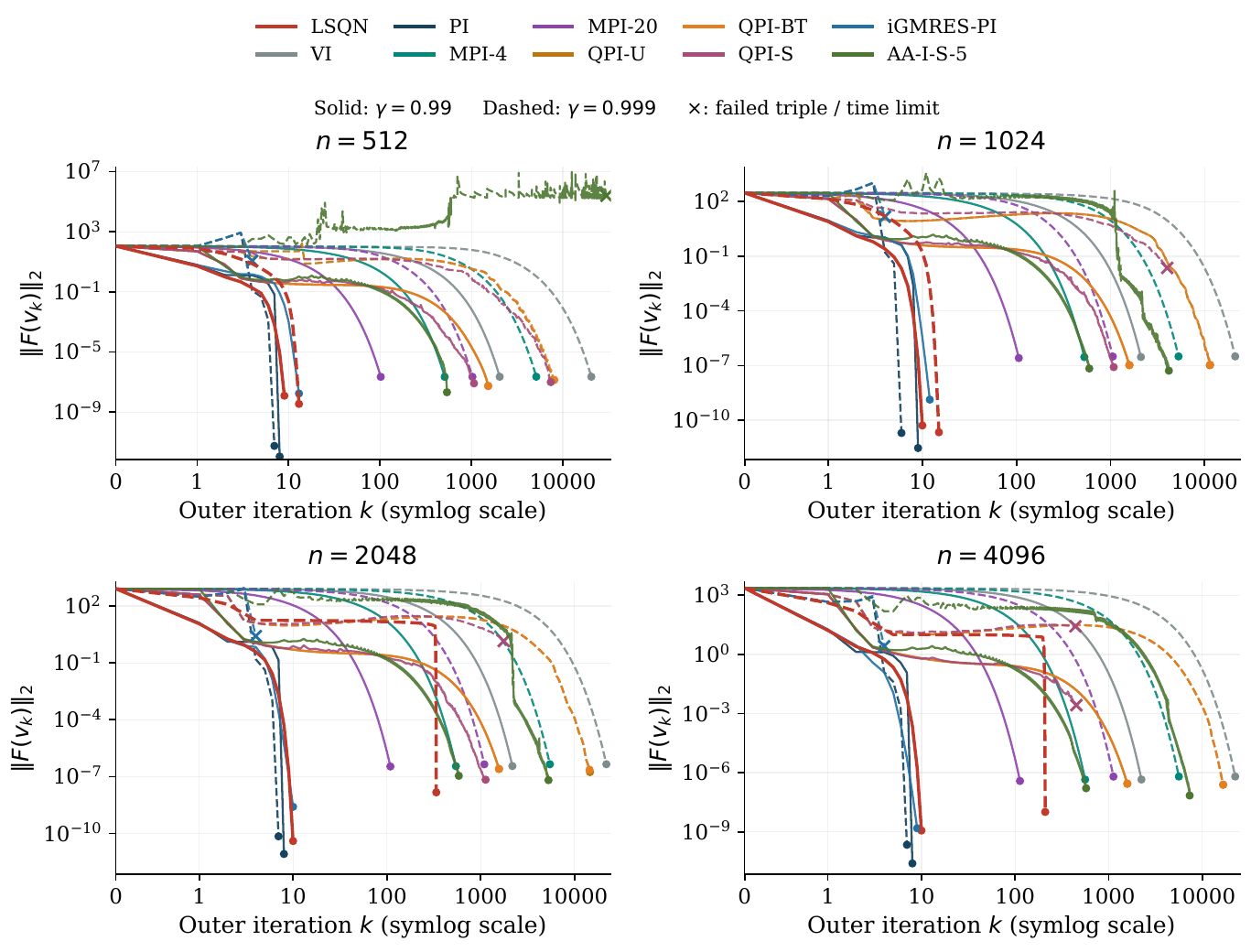}
        \caption{Controlled Queue}\label{fig:controlled_queue}
    \end{minipage}
    \hfill
    \begin{minipage}[b]{0.48\linewidth}
        \centering
        \includegraphics[width=\linewidth]{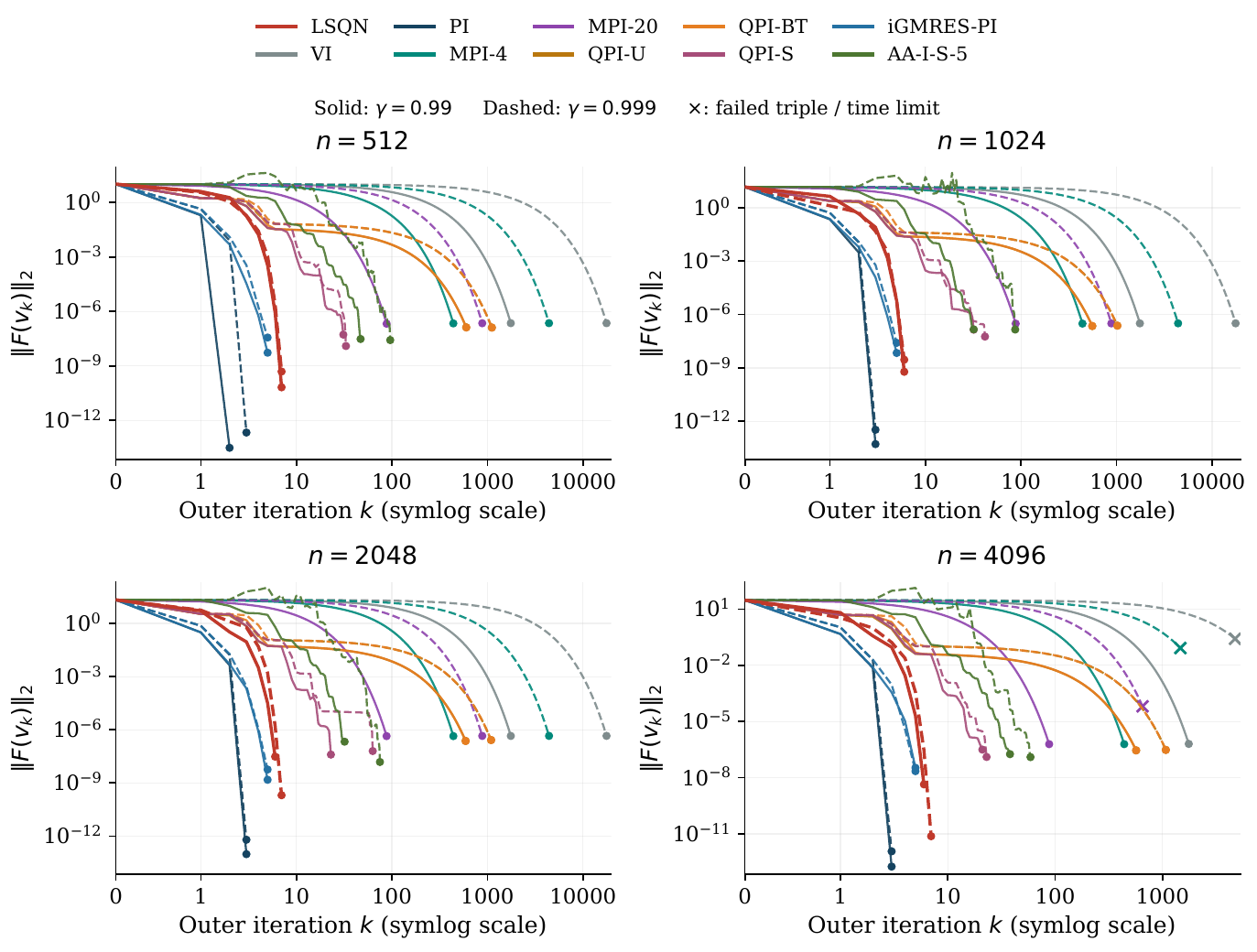}
        \caption{Clustered Chain}\label{fig:clustered_chain}
    \end{minipage}

    \vspace{1em}

    \begin{minipage}[b]{0.48\linewidth}
        \centering
        \includegraphics[width=\linewidth]{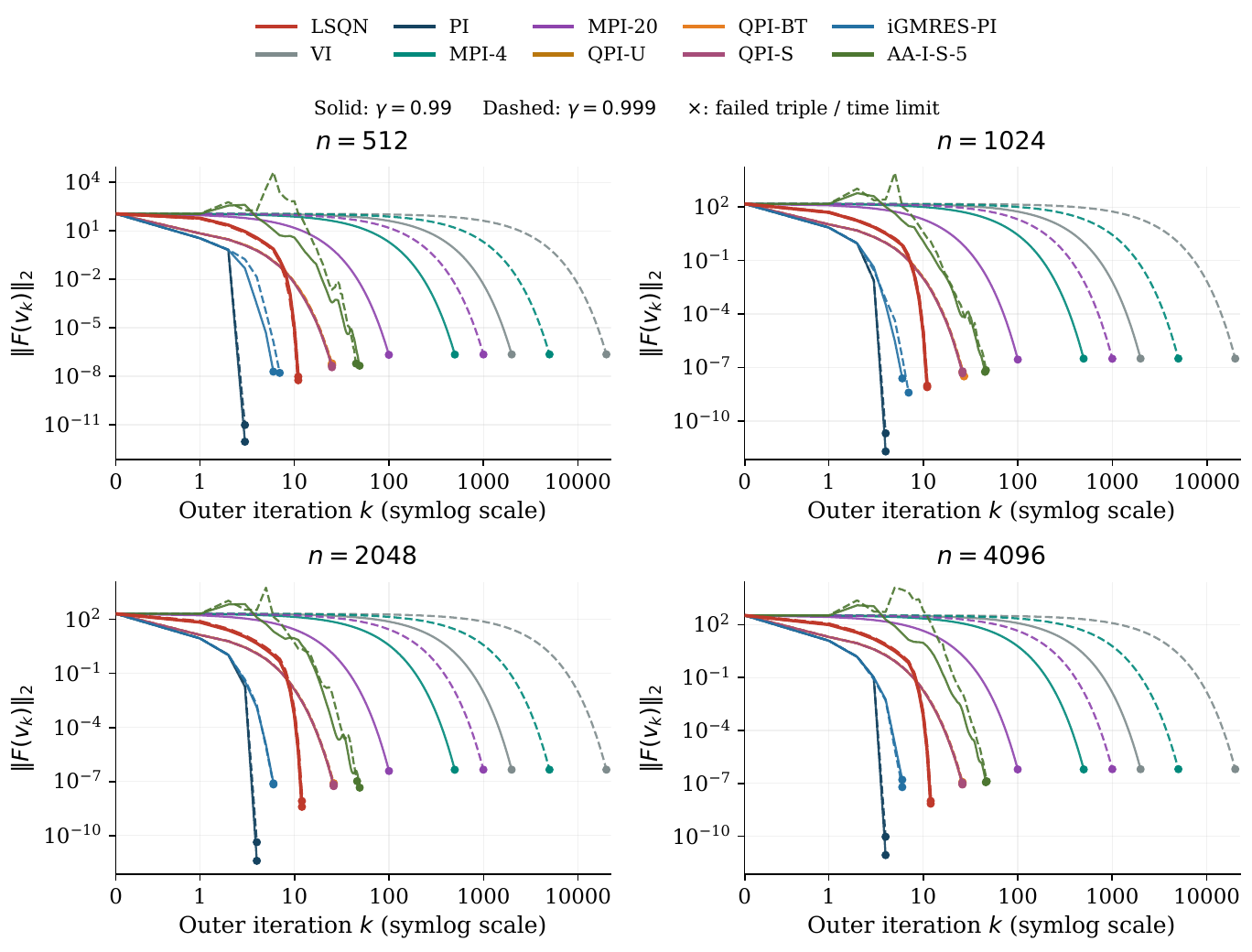}
        \caption{Sparse Garnet}\label{fig:sparse_garnet}
    \end{minipage}
    \hfill
    \begin{minipage}[b]{0.48\linewidth}
        \centering
        \includegraphics[width=\linewidth]{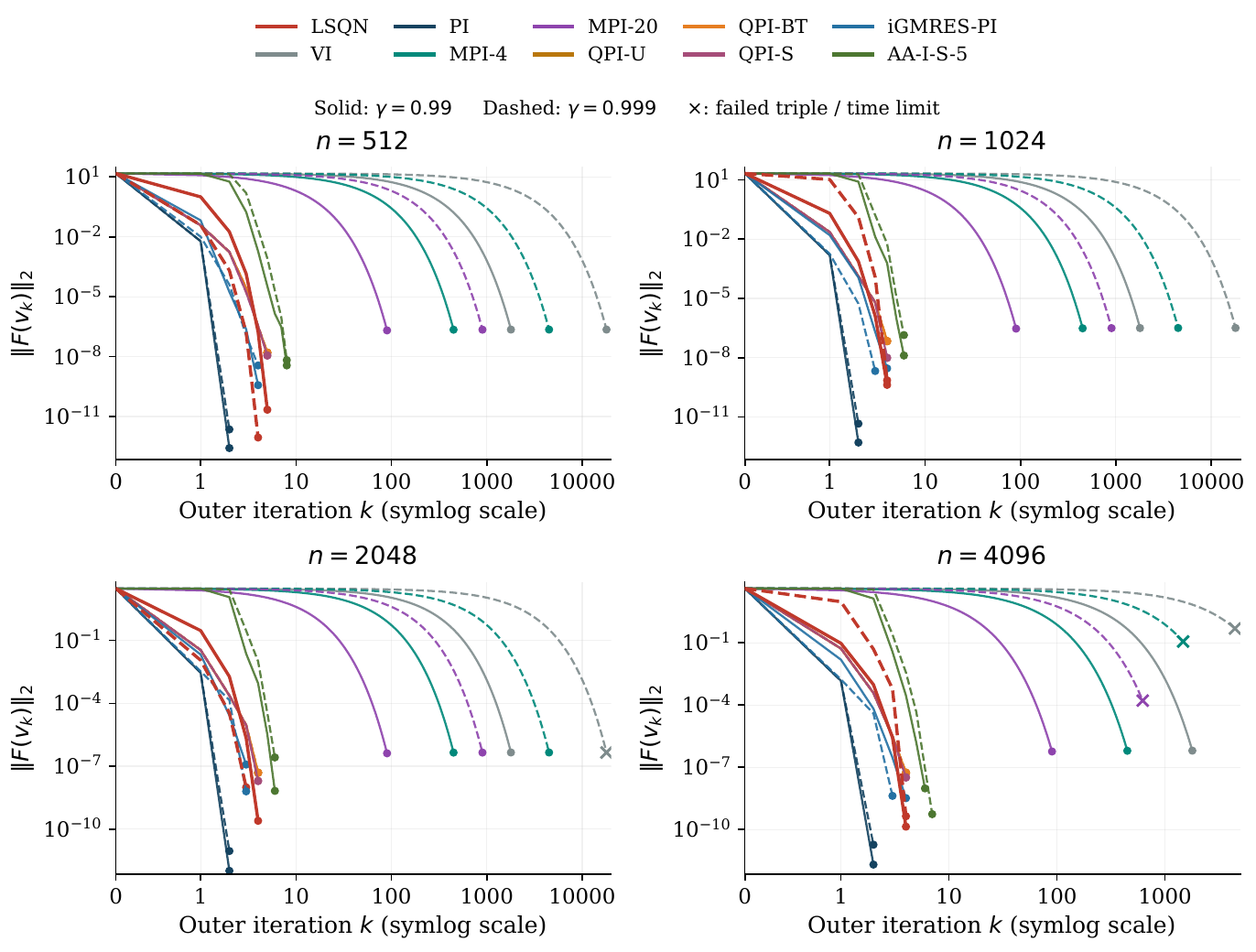}
        \caption{Dense Latent-mixture}\label{fig:dense_latent}
    \end{minipage}
    \label{fig:four_figures}
\end{figure}

We first illustrate the results for the number of outer iterations in testing the $\ell_2$ norm of $F(v_k)$. We present the curve induced by $\gamma = 0.99$ and $\gamma = 0.999$ by solid and dashed lines, respectively in each figure. As shown in Figures \ref{fig:controlled_queue}--\ref{fig:dense_latent}, our algorithm requires generally fewer outer iterations than the other baselines, except for PI and iGMRES-PI, which compute exact and inexact Newton steps, respectively, by solving the current policy equation. 
We also find that as the state number $n$ increases from 512 to 4096 with $\gamma = 0.99$ or $\gamma = 0.999$, the number of outer iterations remains relatively stable. Moreover, the observed results are consistent with the global convergence in Theorem~\ref{thm:global} and the local superlinear convergence in Theorem~\ref{thm:local}.

During the experiment, we add extra checks for failed runs and time limits. Specifically, we attempt three runs for each instance. Only all runs successfully pass the terminal condition is considered to be a feasible instance result. For time limits, we set at most 30 and 60 seconds for $n \le1024$ and $n > 1024$, respectively. For instance, as shown in Figure \ref{fig:controlled_queue}, the QPI-S instance under $n=4096$ and $\gamma=0.999$ does not meet the time limit of 60 seconds in all three runs, thus we give it a dashed line with a cross marker.

Then we focus on runtime and test all 32 instances with each instance given three runs and selecting the median time as runtime. From Table \ref{tab:runtime-scb-final}, LSQN reaches the threshold target in 96 runs and is the fastest on a quarter of the 32 instances, mainly illustrating advantages in the clustered chain model. The algorithm also ranks third and fourth on Sparse Garnet and Dense Latent-mixture, respectively. For PI and its variant, LSQN is faster than exact PI on 24 of the 32 instances, while LSQN is faster than iGMRES-PI on 25 of 28 instances among instances successfully solved by both methods. Moreover, for QPI-S, which also achieves global convergence and local superlinear convergence, LSQN is faster in all 28 instances successfully solved by both methods and also solves the other four instances on which QPI-S reaches the time limit. However, the tridiagonal policy-evaluation systems structure of the Controlled Queue limits LSQN's speed when maintaining a dense inverse model, making inverse-model products and updates costly on the larger high-discount instances. In contrast, note that PI performs well in this model since it uses sparse direct solves. In Table \ref{tab:runtime-scb-final}, successful methods are ranked by median runtime, whereas unsuccessful methods share the average of the remaining ranks.
\begin{table}[t]
\centering
\begingroup
\fontsize{8}{9}\selectfont
\setlength{\tabcolsep}{3.8pt}
\renewcommand{\arraystretch}{1.08}
\begin{tabular}{@{}lrr*{10}{r}@{}}
\toprule
Problem & $n$ & $\gamma$ & \textbf{LSQN} & VI & PI & MPI-4 & MPI-20 & QPI-U & QPI-BT & QPI-S & iGMRES-PI & AA-I-S-5 \\
\midrule
Queue & 512 & 0.99 & 0.0686 & 0.0913 & \textbf{0.0095} & 0.3627 & 0.0819 & 0.1231 & 0.1083 & 1.6564 & 0.1334 & 0.0557 \\
 &  & 0.999 & 0.0807 & 0.9008 & \textbf{0.0070} & 3.7084 & 0.8271 & 0.6086 & 1.9123 & 12.2657 & TO & F(0/3) \\
 & 1024 & 0.99 & 0.2730 & 0.1235 & \textbf{0.0108} & 0.3938 & 0.0924 & 0.1570 & 0.1416 & 8.1624 & 0.1648 & 0.0811 \\
 &  & 0.999 & 0.4710 & 1.2390 & \textbf{0.0073} & 4.0788 & 0.9317 & 1.1310 & 2.6634 & TO & TO & 0.5668 \\
 & 2048 & 0.99 & 0.5219 & 0.2096 & \textbf{0.0131} & 0.4503 & 0.1127 & 0.2128 & 0.1950 & 35.9137 & 0.1981 & 0.1866 \\
 &  & 0.999 & 6.8800 & 1.9889 & \textbf{0.0113} & 4.6928 & 1.1415 & 2.0538 & 3.0101 & TO & TO & 1.2727 \\
 & 4096 & 0.99 & 0.8099 & 0.3668 & \textbf{0.0177} & 0.5883 & 0.1448 & 0.3720 & 0.3487 & TO & 0.3105 & 0.2391 \\
 &  & 0.999 & 26.9971 & 3.6819 & \textbf{0.0154} & 5.9282 & 1.4268 & 3.9456 & 3.7362 & TO & TO & 3.0200 \\
\midrule
Cluster & 512 & 0.99 & \textbf{0.0084} & 0.3660 & 0.0086 & 0.2656 & 0.1229 & 0.1465 & 0.1420 & 0.0583 & 0.0114 & 0.0141 \\
 &  & 0.999 & \textbf{0.0090} & 3.6804 & 0.0132 & 2.6734 & 1.2522 & 0.2718 & 0.2616 & 0.0623 & 0.0143 & 0.0285 \\
 & 1024 & 0.99 & \textbf{0.0255} & 1.4226 & 0.0775 & 1.4510 & 0.5529 & 0.4859 & 0.4806 & 0.2918 & 0.0454 & 0.0339 \\
 &  & 0.999 & \textbf{0.0294} & 14.1915 & 0.0738 & 14.1705 & 5.4892 & 0.8586 & 0.8876 & 0.4003 & 0.0514 & 0.0885 \\
 & 2048 & 0.99 & \textbf{0.0903} & 5.5261 & 0.5014 & 4.6893 & 2.0326 & 1.9925 & 1.9293 & 0.8835 & 0.1715 & 0.1198 \\
 &  & 0.999 & \textbf{0.1059} & 55.6762 & 0.5217 & 46.4943 & 19.9631 & 3.5692 & 3.5188 & 3.0827 & 0.1999 & 0.2670 \\
 & 4096 & 0.99 & \textbf{0.3225} & 21.4234 & 4.0320 & 19.7953 & 8.4453 & 6.9703 & 6.9987 & 4.0106 & 1.2034 & 0.5047 \\
 &  & 0.999 & \textbf{0.4126} & TO & 3.7347 & TO & TO & 13.8129 & 13.5234 & 4.6227 & 1.0582 & 0.7782 \\
\midrule
Garnet & 512 & 0.99 & 0.0030 & 0.1031 & 0.0282 & 0.3897 & 0.0889 & 0.0026 & \textbf{0.0023} & 0.0367 & 0.0264 & 0.0061 \\
 &  & 0.999 & 0.0032 & 1.0393 & 0.0279 & 3.9141 & 0.8929 & 0.0025 & \textbf{0.0022} & 0.0364 & 0.0126 & 0.0050 \\
 & 1024 & 0.99 & 0.0040 & 0.1460 & 0.2276 & 0.4509 & 0.1045 & 0.0033 & \textbf{0.0029} & 0.1911 & 0.0120 & 0.0075 \\
 &  & 0.999 & 0.0043 & 1.4631 & 0.2280 & 4.4497 & 1.0279 & 0.0033 & \textbf{0.0030} & 0.1915 & 0.0138 & 0.0070 \\
 & 2048 & 0.99 & 0.0060 & 0.2643 & 1.7037 & 0.5586 & 0.1356 & \textbf{0.0046} & 0.0048 & 0.8312 & 0.0157 & 0.0164 \\
 &  & 0.999 & 0.0066 & 2.6097 & 1.7375 & 5.4970 & 1.3149 & 0.0051 & \textbf{0.0044} & 0.8387 & 0.0168 & 0.0149 \\
 & 4096 & 0.99 & 0.0101 & 0.4474 & 13.2972 & 0.7759 & 0.1887 & 0.0078 & \textbf{0.0073} & 3.4203 & 0.0180 & 0.0214 \\
 &  & 0.999 & 0.0121 & 4.3800 & 13.3319 & 7.6461 & 1.8568 & 0.0084 & \textbf{0.0072} & 3.4289 & 0.0198 & 0.0224 \\
\midrule
Latent & 512 & 0.99 & 0.0031 & 0.3763 & 0.0084 & 0.2760 & 0.1296 & 0.0017 & \textbf{0.0017} & 0.0089 & 0.0064 & 0.0025 \\
 &  & 0.999 & 0.0028 & 3.8107 & 0.0104 & 2.6855 & 1.2602 & 0.0017 & \textbf{0.0016} & 0.0106 & 0.0067 & 0.0024 \\
 & 1024 & 0.99 & 0.0079 & 1.4637 & 0.0546 & 1.4426 & 0.5730 & 0.0044 & \textbf{0.0044} & 0.0378 & 0.0264 & 0.0065 \\
 &  & 0.999 & 0.0088 & 14.3604 & 0.0554 & 13.9815 & 5.4525 & \textbf{0.0044} & 0.0044 & 0.0359 & 0.0213 & 0.0073 \\
 & 2048 & 0.99 & 0.0314 & 5.5946 & 0.3667 & 4.7967 & 2.0799 & 0.0163 & \textbf{0.0159} & 0.1491 & 0.0770 & 0.0239 \\
 &  & 0.999 & 0.0288 & F(2/3) & 0.3658 & 50.9477 & 21.4017 & 0.0170 & \textbf{0.0167} & 0.1567 & 0.0730 & 0.0231 \\
 & 4096 & 0.99 & 0.1404 & 21.6557 & 2.5239 & 17.8190 & 7.9514 & \textbf{0.0645} & 0.0666 & 0.7618 & 0.4790 & 0.0946 \\
 &  & 0.999 & 0.1105 & TO & 2.9200 & TO & TO & 0.0636 & \textbf{0.0626} & 0.8865 & 0.3115 & 0.1071 \\
\midrule
\multicolumn{3}{@{}l}{Mean rank $\downarrow$} & 3.5312 & 8.1875 & 4.8750 & 8.8438 & 6.4375 & 4.0625 & 3.5625 & 7.0156 & 5.0312 & 3.4531 \\
\multicolumn{3}{@{}l}{Overall rank $\downarrow$} & 2 & 9 & 5 & 10 & 7 & 4 & 3 & 8 & 6 & 1 \\
\multicolumn{3}{@{}l}{Successful instances} & 32/32 & 29/32 & 32/32 & 30/32 & 30/32 & 32/32 & 32/32 & 28/32 & 28/32 & 31/32 \\
\bottomrule
\end{tabular}
\endgroup
\caption{Runtime comparison on 32 discounted MDP instances (seconds).}
\label{tab:runtime-scb-final}
\vspace{2pt}
\begin{minipage}{0.96\linewidth}
\fontsize{8}{9}\selectfont
\end{minipage}
\end{table}

Overall, we compare all 10 algorithms by their speed rankings on the four MDP models and make an average. Our algorithm ranks second with a mean rank of 3.53, showing a huge potential for solving MDPs. Although basic algorithm QPI-U and its variant QPI-BT perform well on the Sparse Garnet and Dense Latent-mixture models, LSQN represents the strongest performance especially on the Clustered Chain model. In addition, as a quasi-Newton based method, we have validated its potential for accelerating VI compared with MPI-4 and MPI-20 which can be viewed as methods between VI and PI. Note that in the Table \ref{tab:runtime-scb-final}, "TO" means exceeding time limits under this method, while "F($s$/3)" means this scheme runs $s$ times successfully through all three runs. 

\subsection{Inner Learning on Public Benchmarks}
\label{subsec:public-inner-learning}

We next examine the number of learning queries required by LSQN relative
to the exact-arithmetic worst-case bound of $n$ queries per outer iteration.
Recall that $m_k$, defined in Subsection~\ref{subsec:wccc}, denotes the number of inner queries at iteration $k$. 
In this subsection, we demonstrate empirically that $m_k$ is typically much smaller than $n$.

We first consider the FrozenLake, Taxi, and CliffWalking benchmark
environments provided by Gymnasium \citep{Towers2024Gymnasium}.
The transition models are extracted from Gymnasium~1.1.1, using the
default slippery $4\times4$ and $8\times8$ FrozenLake maps, Taxi-v3,
and CliffWalking-v0. An additional zero-reward absorbing state represents
episode termination while preserving the reward on terminal entry.
The designed state counts are $17$, $65$, $501$, and $49$, respectively.
Taxi has six actions and the other environments have four.
We use $\gamma\in\{0.99,0.999\}$ as well, and the initial inverse model by $H_0=I$.
As in subsection \ref{subsec:runtime}, each instance is run three times with a 30-second time limit. Success is checked independently using
$\|F(v)\|_\infty\leq10^{-8}$.

\begin{table}[htbp]
\centering
\small
\setlength{\tabcolsep}{5pt}
\renewcommand{\arraystretch}{1.08}
\begin{tabular*}{\linewidth}{@{\extracolsep{\fill}}lrrrrrrr@{}}
\toprule
& & \multicolumn{3}{c}{$\gamma=0.99$}
  & \multicolumn{3}{c}{$\gamma=0.999$} \\
\cmidrule(lr){3-5}\cmidrule(l){6-8}
Environment & $n$
& $\bar m$ & $m_{\max}$ & Time (s)
& $\bar m$ & $m_{\max}$ & Time (s) \\
\midrule
FrozenLake $4\times4$ & 17
& 2.71 & 5 & 0.0066 & 2.71 & 5 & 0.0059 \\
FrozenLake $8\times8$ & 65
& 5.25 & 9 & 0.0113 & 5.43 & 10 & 0.0135 \\
Taxi & 501
& 24.17 & 100 & 0.2364 & 75.94 & 353 & 0.9054 \\
CliffWalking & 49
& 1.76 & 4 & 0.0143 & 2.96 & 6 & 0.0201 \\
\bottomrule
\end{tabular*}
\caption{Learning queries on the official benchmarks.
$\bar m$ and $m_{\max}$ are the mean and maximum number of inner queries per outer iteration. 
 Time denotes the median runtime of the complete solver over three runs.}
\label{tab:public-learning-queries}
\end{table}

All 24 runs complete successfully. Counting each configuration once, LSQN completes 174 outer iterations, with none doing inner learning using all $n$ learning queries. Moreover, in 13 iterations, the proposal generated by the incoming inverse model already satisfies the forcing test, requiring no additional learning queries. Table~\ref{tab:public-learning-queries} shows that the benchmark FrozenLake and CliffWalking require only a few queries on average. Instead, Taxi is a little more demanding, particularly at $\gamma=0.999$, where the mean is $75.94$
and the maximum is $353$, but still below $n=501$.
The proposed results show that the forcing test can often be passed before full-dimensional learning queries. 

Additionally, in order to illustrate growth with dimension within one model, we also report a separate synthetic FrozenLake scaling experiment. These models
are local reconstructions of the transition rules, rather than additional runs of the official environments above. 
The sizes are considered to be $n\in\{64,256,576,1024,2304,4096\}$, with three random seeds 130, 131 and 132 at each state size. For a fixed discount factor, let $K_{n,s}$ denote the number of outer iterations
for size $n$ and seed $s$. We first average queries within each instance
and then average the three instance means with equal weights by the rule: $\bar m(n)=\frac{1}{3}\sum_{s\in\{130,131,132\}}
\left(\frac{1}{K_{n,s}}\sum_{k=0}^{K_{n,s}-1}m_k^{(n,s)}\right).$

\begin{table}[htbp]
\centering
\small
\setlength{\tabcolsep}{5pt}
\renewcommand{\arraystretch}{1.08}
\begin{tabular*}{\linewidth}{@{\extracolsep{\fill}}lrrrrrr@{}}
\toprule
$\gamma$ & $n=64$ & $n=256$ & $n=576$
& $n=1024$ & $n=2304$ & $n=4096$ \\
\midrule
0.99  & 4.92 & 6.72 & 8.08  & 8.78  & 9.68  & 9.42 \\
0.999 & 5.23 & 8.46 & 10.63 & 14.94 & 17.76 & 20.22 \\
\bottomrule
\end{tabular*}
\caption{Mean learning-query counts $\bar m(n)$ for synthetic FrozenLake
maps. Each entry averages three instance means with equal weights.}
\label{tab:learning-scaling}
\end{table}

We then fit $\log\bar m(n)=a+\rho\log n$ by ordinary least squares over
the six state sizes, separately for each discount and using unrounded means.
The resulting empirical power laws are
\begin{equation}
\bar m(n)=
\begin{cases}
2.654\,n^{0.1643}, & \gamma=0.99,\\
1.325\,n^{0.3342}, & \gamma=0.999.
\end{cases}
\label{eq:empirical-learning-growth}
\end{equation}
Moreover, Table \ref{tab:learning-scaling} shows that the average number of learning queries grows much more slowly than the problem dimension \(n\). Specifically, as \(n\) increases from \(64\) to \(4096\), the mean query count rises from \(4.92\) to \(9.42\) for \(\gamma=0.99\), and from \(5.23\) to \(20.22\) for \(\gamma=0.999\). 
In particular, based on the empirical fit in \eqref{eq:empirical-learning-growth}, the estimated average inner-learning computational cost grows approximately as \(n^{2.1643}\) and \(n^{2.3342}\) over the tested range of $n$ for the two discount factors, respectively.

\FloatBarrier
\section{Conclusion}
This paper establishes a quasi-Newton method for MDPs, where an Armijo line-search scheme is applied to guarantee the monotonic decrease of the $\ell_2$ norm of the Bellman residual without requiring the uniqueness of the optimal policy, achieving global linear convergence and local superlinear convergence rate. Numerical experiments are conducted to evaluate the convergence behavior and illustrate the advantages of this method to implement in solving MDPs, which implies its potential to accelerate VI.

\appendix
\section{Clarke Containment and Envelope Support}\label{app:clarke}
This appendix derives the containment and support properties of the generalized-gradient envelope used in LSQN, using standard Clarke calculus \citep{Clarke1990}.

Write $\ell_i^a(v)=v_i-r_i^a-\gamma P_i^av$ and let $F_i(v)$ be the $i$-th component of $F(v)$. By \eqref{eq:q_definitions},
\begin{equation}
 F_i(v)=\min_{a\in\A}\ell_i^a(v),\qquad
 \nabla\ell_i^a=e_i-\gamma(P_i^a)^\top,
 \label{eq:minimum-residual}
\end{equation}
where $e_i\in\mathbb R^n$ denotes the $i$-th standard basis
vector. Let $\A_i^0(v)=\argmax_a q_i^a(v)$ and by finiteness of the action set $\A_i^0(v) \neq \emptyset$.
It follows from (\ref{eq:q_definitions}) that the elements of $\A_i^0(v)$ minimize $\ell_i^a(v)$ and also belong to $\A_i^\varepsilon(v)$.
Applying the finite-maximum rule to $-F_i$ and then negating gives
\begin{equation}
 \partial_CF_i(v)
 =-\partial_C\max_a\{-\ell_i^a(v)\}=-\operatorname{conv}
       \{-\nabla\ell_i^a:a\in\A_i^0(v)\}=\operatorname{conv}
       \{e_i-\gamma(P_i^a)^\top:a\in\A_i^0(v)\}.
 \label{eq:component-clarke}
\end{equation}
For $\phi_i(v)=F_i(v)^2/2$, the Clarke chain rule yields
$\partial_C\phi_i(v)\subseteq F_i(v)\partial_CF_i(v)$.
The sum rule and \eqref{eq:component-clarke} therefore imply
 \begin{align}
 \partial_C\Phi(v)
 &\subseteq\sum_i F_i(v)\partial_CF_i(v)=\left\{\sum_i F_i(v)\sum_{a\in\A_i^0(v)}
             \lambda_{i}^a\bigl(e_i-\gamma(P_i^a)^\top\bigr):
                  \lambda_i\in\Delta(\A_i^0(v))\right\}  \label{eq:clarke-containment}\\
 &=\{J_\lambda^\top F(v):\lambda_i\in\Delta(\A_i^0(v))\}\subseteq\G_\varepsilon(v),\notag
\end{align}
where $\sum_i F_i(v)\partial_CF_i(v)$ is a Minkowski sum and the final inclusion uses the fact that
$\A_i^0(v)\subseteq\A_i^\varepsilon(v)$.

The product $\prod_i\Delta(\A_i^\varepsilon(v))$ is nonempty, compact and
convex because every near-active set is finite and contains a greedy action.
Its image under the affine map $\lambda\mapsto J_\lambda^\top F(v)$ is
therefore nonempty, compact and convex. For any $d \in \mathbb{R}^n$,
\begin{align}
&\max_{h\in\widehat{\mathcal G}_\varepsilon(v)} h^\top d=
\max_{\substack{
    \lambda_i\in\Delta(\mathcal A_i^\varepsilon(v))\\
    i=1,\ldots,n}}
\sum_{i=1}^n
\sum_{a\in\mathcal A_i^\varepsilon(v)}
\lambda_i^a F_i(v)
\bigl(d_i-\gamma P_i^a d\bigr) \notag
\\&=
\sum_{i=1}^n
\max_{\lambda_i\in\Delta(\mathcal A_i^\varepsilon(v))}
\sum_{a\in\mathcal A_i^\varepsilon(v)}
\lambda_i^a F_i(v)
\bigl(d_i-\gamma P_i^a d\bigr)
=
\sum_{i=1}^n
\max_{a\in\mathcal A_i^\varepsilon(v)}
F_i(v)\bigl(d_i-\gamma P_i^a d\bigr),\label{eq:separable-support-proof}
\end{align}
where the first equality in \eqref{eq:separable-support-proof} uses independence across rows, and the second uses linearity over a simplex. Replacing maxima by minima gives \eqref{eq:M}
by the same arguments. Finally, the support representation of the Clarke
directional derivative and \eqref{eq:clarke-containment} give
\begin{equation}
 \Phi^\circ(v;d)
 =\max_{h\in\partial_C\Phi(v)}h^\top d
 \le\max_{h\in\G_\varepsilon(v)}h^\top d
 =D(v;d).
 \label{eq:clarke-directional}
\end{equation}
It is noteworthy that the set $\G_{\varepsilon}(v)$ in equation (\ref{eq:clarke-containment}) is constructed by independently selecting the active weights $\lambda_i$ for each row $i$ and combining them. However, in practice, we can only optimize a unified policy through a single perturbation direction. This means it is difficult to guarantee that a single perturbation can simultaneously \textit{jointly realize} the ideal greedy policy for every row.

Because of this potential conflict, we can only establish a containment relationship (i.e., $\partial_C\Phi(v) \subseteq \G_{\varepsilon}(v)$). Nevertheless, this containment is entirely sufficient for the subsequent derivation: the global step-length argument in Lemma~\ref{lem:capture} is obtained directly along specific segments and does not rely on the strict equality of these two sets.

\section{Test Models and Benchmark Protocol}
\label{app:models}

This appendix specifies the models, solver settings, and reporting conventions
used in Section \ref{sec:experiments}, including the four-family comparison, and the official public benchmarks experiment. Throughout the appendix, $\mathbf{1}_m = (1,...,1)^\top\in\mathbb{R}^m$
denotes the all-ones column vector, and $\mathbb{I}_{\{\cdot \}}$ is the indicator.

\subsection{Synthetic Test Models}
\label{app:synthetic-models}

The main comparison contains four model families, each with four actions,
$n\in\{512,1024,2048,4096\}$, and
$\gamma\in\{0.99,0.999\}$, giving 32 instances.
Each model is generated with NumPy's random generator initialized with seed 130.
For a fixed family and state count, the transition matrices and rewards are
identical for the two discount factors. The descriptions below use zero-based
state indices $i,j\in\{0,\ldots,n-1\}$ and action indices
$a\in\{0,1,2,3\}$. All Gaussian perturbations $\omega_i^a,\psi_{i\ell}, \varphi_{i\ell}^a,\upsilon_i^a \overset{\mathrm{iid}}{\sim}\mathcal{N}(0,1)$ are independent standard normal
draws, except where a shared array is explicitly specified.

\paragraph{Controlled queue.}
States represent queue lengths. The arrival probability is $0.30$, and the
service probability under action $a$ is $s_a=0.08+0.18a$.
Set $u_a=0.30(1-s_a)$, $d_a=0.70s_a$, and $h_a=1-u_a-d_a$.
The transition probabilities are defined as
\[
 P_{ij}^a
 =u_a\mathbb{I}_{\{j=\min(n-1,i+1)\}}
 +d_a\mathbb{I}_{\{j=\max(0,i-1)\}}
 +h_a\mathbb{I}_{\{j=i\}}.
\]
Thus, transitions leaving the state space are accumulated at the corresponding
boundary. Moreover, rewards are
\(r_i^a=-0.015i-0.11a^{1.4}+0.002\omega_i^a.\)

\paragraph{Nearly decomposable clustered chain.}
Partition the states into six consecutive groups using
$c_i=\lfloor6i/n\rfloor$, and define the group as an ordered sequence \(C_k = (j_1, j_2, \dots, j_{|C_k|})\), where \(j_1 < j_2 < \dots < j_{|C_k|}\) are all the states \(j\) such that \(c_j = k\). For each pair of states and actions, we draw a probability vector
$w_i^a\sim\operatorname{Dirichlet}(0.4\mathbf{1}_{|C_{c_i}|})$ on $C_{c_i}$, extending it by zero outside the group, and sample a set $S_i^a$
of three distinct states uniformly outside $C_{c_i}$.
If $i$ has zero-based position $t_i$ in the ordered group $C_{c_i}$, let $d_i^a$ be
the group member at position $(t_i+a+1)\bmod |C_{c_i}|$. Define the transition probability as
\[
 P_{ij}^a=0.92\widetilde P_{ij}^a
             +0.08\mathbb{I}_{\{j=d_i^a\}}, \quad
 \widetilde P_{ij}^a
 =0.99w_{ij}^a\mathbb{I}_{\{j\in C_{c_i}\}}
  +\frac{0.01}{3}\mathbb{I}_{\{j\in S_i^a\}}.
\]
With $\phi_i=2\pi i/n$, rewards are denoted as
\( r_i^a=-0.05c_i+0.15\cos(\phi_i+a\pi/3)-0.025a+0.01\omega_i^a.
\)
\paragraph{Sparse Garnet model.}
For each pair of states and actions, we select eight distinct successor states without replacement and draw their probabilities from
$\operatorname{Dirichlet}(\mathbf{1}_8)$. All other transition probabilities
are set to be zero, and rewards $r_i^a$ are independent standard normal draws.

\paragraph{Dense latent-mixture model.}
Draw five shared probability rows
$V_{\ell}\sim\operatorname{Dirichlet}(0.35\mathbf{1}_n)$,
$\ell \in \{0,\ldots,4\}$. Using $\phi_i=2\pi i/n$, we construct the shared feature array and action-dependent logits by
\[
 B_{i\ell}=\sin((\ell+1)\phi_i)+0.2\psi_{i\ell},
 \qquad
 Z_{i\ell}^a=1.3B_{i,(\ell-a)\bmod5}+0.25\varphi_{i\ell}^a.
\]
Given \(
 W_{i\ell}^a=
 \frac{\exp(Z_{i\ell}^a)}{\sum_{h=0}^{4}\exp(Z_{ih}^a)}
 \), we obtain its transition probability matrix \(P^a=W^aV,
\)
and use rewards \( r_i^a=0.28\cos(\phi_i-2\pi a/4)+0.08\sin(3\phi_i+a)+0.03\upsilon_i^a.\)

\subsection{Solver Configurations}
\label{app:solver-configurations}

For the main comparison and the official benchmarks, all methods use the same initial value as
\(
 v_0=\frac{r_{\min}}{1-\gamma}\mathbf{1}_n\), where \(r_{\min}=\min_{i,a}r_i^a.
\) The nine baselines are VI, PI, MPI-4, MPI-20, QPI-U, QPI-BT, QPI-S, iGMRES-PI, and AA-I-S-5. The settings of LSQN have been mentioned in subsection \ref{subsec:experiment_setup}.

\paragraph{VI, PI, and MPI.}
VI simply applies the Bellman operator, while PI evaluates the selected policy using a dense
linear solve or a sparse direct solve according to the transition storage.
MPI-4 and MPI-20 use four and twenty policy-evaluation sweeps, respectively. Note that for MPI, the policy matrix is assembled once per outer iteration and reused
throughout these sweeps.

\paragraph{QPI variants.}
Our implementations follow the equations of
\citet{KolarijaniEsfahani2023QPI}.
QPI-U uses a uniform transition prior and the original residual safeguard,
with threshold $\gamma^{k+1}\|F(v_0)\|_\infty$.
QPI-BT uses a trial path anchored at $Tv_k$, a backtracking factor of $1/2$,
and a residual reduction factor of $(1+\gamma)/2$.
QPI-S augments this construction with the policy-dependent secant constraints,
using a least-Frobenius-change transition update and a Woodbury inverse update.
Its constraint calculation uses a relative singular-value cutoff of $10^{-12}$.
An update is reported as unsuccessful if the constraint residual exceeds
$10^{-8}\max\{1,\|B_{\mathrm{con}}\|_F\}$, where $B_{\mathrm{con}}$ is the
constraint right-hand side, or if the small Woodbury system has condition
number greater than $10^{14}$. Numerical failures are retained in the results.

\paragraph{iGMRES-PI.}
We execute the authors' \texttt{GMRES\_PI} Python routine associated with
\citet{Gargiani2024IPI}, using forcing parameter $0.1$ and GMRES restart
length 100. The adapter converts rewards to costs, supplies the common
dense or sparse policy matrices, and uses the common Bellman-residual stopping
test. The authors' tolerance conversion is retained: the relative Euclidean
GMRES tolerance is $0.1$ times the initial policy residual in $\ell_{\infty}$ norm, divided by the Euclidean norm of the policy reward vector.

\paragraph{AA-I-S-5.}
This baseline is a Python port of the authors' stabilized type-I Anderson
implementation \citep{ZhangODonoghueBoyd2020}, with memory $m=5$, Powell
regularization, restart, and the original declining-threshold safeguard.
Its parameters are
$\theta_{\mathrm{AA}}=0.01$, $\tau_{\mathrm{AA}}=0.001$,
$D_{\mathrm{AA}}=10^6$, and $\varepsilon_{\mathrm{AA}}=10^{-6}$.
The fallback is the damped VI step conducted by $0.1v+0.9Tv$.

\subsection{Public Benchmarks and Synthetic Scaling}
\label{app:public-scaling-protocol}

\paragraph{Official public environments.}
We extract the transition tables directly from Gymnasium~1.1.1
\citep{Towers2024Gymnasium}, using the default slippery FrozenLake
$4\times4$ and $8\times8$ maps, Taxi-v3, and CliffWalking-v0.
Every terminal transition is redirected to one additional zero-reward
absorbing state, while its entry reward is preserved. The resulting state
counts are 17, 65, 501, and 49, while Taxi has six actions and the other environments
have four, as also mentioned in subsection \ref{subsec:public-inner-learning}. 
For each environment and $\gamma\in\{0.99,0.999\}$, LSQN starts from $v_0$ as specified in \ref{app:solver-configurations} among all states, including the added
absorbing state. Each configuration is run three times, giving 24 runs.
Note that reported runtimes are median whole-solver times.

\paragraph{Synthetic scaling experiment.}
For the scaling experiment in subsection \ref{subsec:public-inner-learning}, we generate square FrozenLake maps with the start and goal at opposite corners. Each state other than the start and goal is independently terminal with probability $0.1$. Transitions follow the slippery FrozenLake rules \citep{Towers2024Gymnasium}. Terminal states are absorbing, and rewards are zero except for a reward of one upon entering the goal. Moreover, no additional absorbing state is introduced, and $v_0=0$.
\bibliographystyle{plainnat}
\setlength{\bibsep}{3pt}
\bibliography{references}

@book{Puterman1994,
  author    = {Martin L. Puterman},
  title     = {{Markov} Decision Processes: Discrete Stochastic Dynamic Programming},
  publisher = {Wiley},
  address   = {New York, NY, USA},
  year      = {1994}
}

@book{Bertsekas2012,
  author    = {Dimitri P. Bertsekas},
  title     = {Dynamic Programming and Optimal Control, Vol. II},
  edition   = {4th},
  publisher = {Athena Scientific},
  year      = {2012}
}

@article{Broyden1965,
  author  = {C. G. Broyden},
  title   = {A Class of Methods for Solving Nonlinear Simultaneous Equations},
  journal = {Mathematics of Computation},
  volume  = {19},
  number  = {92},
  pages   = {577--593},
  year    = {1965},
  doi     = {10.1090/S0025-5718-1965-0198670-6}
}

@book{Clarke1990,
  author    = {Frank H. Clarke},
  title     = {Optimization and Nonsmooth Analysis},
  publisher = {SIAM},
  address   = {Philadelphia, PA, USA},
  year      = {1990}
}

@inproceedings{LiuEtAl2023BlockBroyden,
  author    = {Chengchang Liu and Cheng Chen and Luo Luo and John C. S. Lui},
  title     = {Block {B}royden's Methods for Solving Nonlinear Equations},
  booktitle = {Advances in Neural Information Processing Systems},
  volume    = {36},
  year      = {2023}
}

@article{BurdakovKamandi2018,
  author  = {Oleg Burdakov and Ahmad Kamandi},
  title   = {Multipoint Secant and Interpolation Methods with Nonmonotone Line Search for Solving Systems of Nonlinear Equations},
  journal = {Applied Mathematics and Computation},
  volume  = {338},
  pages   = {421--431},
  year    = {2018},
  doi     = {10.1016/j.amc.2018.05.041}
}

@article{Gargiani2024SSN,
  author  = {Matilde Gargiani and Andrea Zanelli and Dominic Liao-McPherson and Tyler Summers and John Lygeros},
  title   = {Dynamic Programming Through the Lens of Semismooth {N}ewton-Type Methods},
  journal = {IEEE Control Systems Letters},
  volume  = {6},
  pages   = {2996--3001},
  year    = {2022},
  note    = {Extended version: arXiv:2203.08678}
}

@article{Gargiani2024IPI,
  author  = {Matilde Gargiani and Robin Sieber and Efe Balta and Dominic Liao-McPherson and John Lygeros},
  title   = {Inexact Policy Iteration Methods for Large-Scale {Markov} Decision Processes},
  journal = {arXiv preprint arXiv:2404.06136},
  year    = {2024}
}

@article{Gargiani2023GMRES,
  author    = {Matilde Gargiani and Dominic Liao-McPherson and Andrea Zanelli and John Lygeros},
  title     = {Inexact {GMRES} Policy Iteration for Large-Scale {Markov} Decision Processes},
  journal   = {IFAC-PapersOnLine},
  volume    = {56},
  number    = {2},
  pages     = {11249--11254},
  year      = {2023},
  doi       = {10.1016/j.ifacol.2023.10.316}
}

@article{GrandClement2021Review,
  author  = {Julien Grand-Cl{\'e}ment},
  title   = {From Convex Optimization to {MDP}s: A Review of First-Order, Second-Order and Quasi-{N}ewton Methods for {MDP}s},
  journal = {arXiv preprint arXiv:2104.10677},
  year    = {2021},
  doi     = {10.48550/arXiv.2104.10677}
}

@article{GoyalGrandClement2022AVI,
  author  = {Vineet Goyal and Julien Grand-Cl{\'e}ment},
  title   = {A First-Order Approach to Accelerated Value Iteration},
  journal = {Operations Research},
  volume  = {71},
  number  = {2},
  pages   = {517--535},
  year    = {2023},
  doi     = {10.1287/opre.2022.2269}
}

@inproceedings{LeeRyu2023AncVI,
  author    = {Jongmin Lee and Ernest K. Ryu},
  title     = {Accelerating Value Iteration with Anchoring},
  booktitle = {Advances in Neural Information Processing Systems},
  volume    = {36},
  year      = {2023}
}

@article{DemboEisenstatSteihaug1982,
  author  = {Ron S. Dembo and Stanley C. Eisenstat and Trond Steihaug},
  title   = {Inexact {Newton} Methods},
  journal = {SIAM Journal on Numerical Analysis},
  volume  = {19},
  number  = {2},
  pages   = {400--408},
  year    = {1982},
  doi     = {10.1137/0719025}
}

@article{ZhangODonoghueBoyd2020,
  author = {Junzi Zhang and Brendan O'Donoghue and Stephen Boyd},
  title = {Globally Convergent Type-{I} {Anderson} Acceleration for Nonsmooth Fixed-Point Iterations},
  journal = {SIAM Journal on Optimization},
  volume = {30}, number = {4}, pages = {3170--3197}, year = {2020},
  doi = {10.1137/18M1232772}
}

@article{ThemelisPatrinos2019,
  author = {Andreas Themelis and Panagiotis Patrinos},
  title = {{SuperMann}: A Superlinearly Convergent Algorithm for Finding Fixed Points of Nonexpansive Operators},
  journal = {IEEE Transactions on Automatic Control},
  volume = {64}, number = {12}, pages = {4875--4890}, year = {2019},
  doi = {10.1109/TAC.2019.2906393}
}

@article{SunHan1997,
  author = {Defeng Sun and Jiye Han},
  title = {{Newton} and Quasi-{Newton} Methods for a Class of Nonsmooth Equations and Related Problems},
  journal = {SIAM Journal on Optimization},
  volume = {7}, number = {2}, pages = {463--480}, year = {1997},
  doi = {10.1137/S1052623494274970}
}

@article{GowerRichtarik2017,
  author = {Robert M. Gower and Peter Richt{\'a}rik},
  title = {Randomized Quasi-{Newton} Updates Are Linearly Convergent Matrix Inversion Algorithms},
  journal = {SIAM Journal on Matrix Analysis and Applications},
  volume = {38}, number = {4}, pages = {1380--1409}, year = {2017},
  doi = {10.1137/16M1062053}
}

@article{FangSaad2009,
  author = {Fang, Haw-ren and Saad, Yousef},
  title = {Two Classes of Multisecant Methods for Nonlinear Acceleration},
  journal = {Numerical Linear Algebra with Applications},
  volume = {16}, number = {3}, pages = {197--221}, year = {2009},
  doi = {10.1002/nla.617}
}

@article{DennisSchnabel1979,
 author={Dennis, Jr., J. E. and Schnabel, Robert B.},
 title={Least Change Secant Updates for Quasi-{Newton} Methods},
 journal={SIAM Review}, volume={21}, number={4}, pages={443--459},
 year={1979}, doi={10.1137/1021091}
}

@article{DennisWalker1985,
 author={Dennis, Jr., J. E. and Walker, Homer F.},
 title={Least-Change Sparse Secant Update Methods with Inaccurate Secant Conditions},
 journal={SIAM Journal on Numerical Analysis}, volume={22}, number={4},
 pages={760--778}, year={1985}, doi={10.1137/0722046}
}

@ARTICLE{11535449,
  author={Zhou, Wangzhi and Mo, Yuanqiu and Dasgupta, Soura},
  journal={IEEE Control Systems Letters}, 
  title={Finite-Time Optimal Policy Identification for the Stochastic Shortest Path Problem}, 
  year={2026},
  volume={10},
  number={},
  pages={409-414},
  doi={10.1109/LCSYS.2026.3697256}}

@article{BroydenDennisMore1973,
  author  = {Broyden, C. G. and Dennis, Jr., J. E. and Mor{\'e}, Jorge J.},
  title   = {On the Local and Superlinear Convergence of Quasi-{Newton} Methods},
  journal = {Journal of the Institute of Mathematics and Its Applications},
  volume  = {12},
  number  = {3},
  pages   = {223--245},
  year    = {1973},
  doi     = {10.1093/imamat/12.3.223}
}

@article{Towers2024Gymnasium,
  author = {Towers, Mark and Kwiatkowski, Ariel and Terry, Jordan
            and Balis, John U. and De Cola, Gianluca and Deleu, Tristan
            and Goul{\~a}o, Manuel and Kallinteris, Andreas
            and Krimmel, Markus and KG, Arjun and Perez-Vicente, Rodrigo
            and Pierr{\'e}, Andrea and Schulhoff, Sander
            and Tai, Jun Jet and Tan, Hannah and Younis, Omar G.},
  title = {{Gymnasium}: A Standard Interface for Reinforcement Learning Environments},
  journal = {arXiv preprint arXiv:2407.17032},
  year = {2024},
  doi = {10.48550/arXiv.2407.17032},
}

@article{KolarijaniEsfahani2023QPI,
  author  = {Sharifi Kolarijani, Mohamad Amin
             and Mohajerin Esfahani, Peyman},
  title   = {From Optimization to Control: Quasi-Policy Iteration},
  journal = {IEEE Transactions on Automatic Control},
  year    = {2026},
  volume  = {71},
  number  = {6},
  pages   = {3880--3893},
  doi     = {10.1109/TAC.2026.3652902}
}

@article{PutermanShin1978,
  author  = {Puterman, Martin L. and Shin, Moon Chirl},
  title   = {Modified Policy Iteration Algorithms for Discounted
             {Markov} Decision Problems},
  journal = {Management Science},
  volume  = {24},
  number  = {11},
  pages   = {1127--1137},
  year    = {1978},
  doi     = {10.1287/mnsc.24.11.1127}
}

@article{KamanchiDiddigiBhatnagar2022,
  author  = {Kamanchi, Chandramouli and
             Diddigi, Raghuram Bharadwaj and
             Bhatnagar, Shalabh},
  title   = {Generalized Second-Order Value Iteration in
             {Markov} Decision Processes},
  journal = {IEEE Transactions on Automatic Control},
  volume  = {67},
  number  = {8},
  pages   = {4241--4247},
  year    = {2022},
  doi     = {10.1109/TAC.2021.3112851}
}
\end{document}